  \documentclass{article}
 \usepackage{amsmath,amsfonts,amssymb,amsthm}
 \usepackage{enumerate,color}
\usepackage[applemac]{inputenc}
 \usepackage[polutonikogreek,english]{babel}

  \date{October 10, 2009}   

\newcommand{\ds}{\displaystyle}

\theoremstyle{plain}
\newtheorem{theorem}{Theorem}[section]
\newtheorem{proposition}[theorem]{Proposition}
\newtheorem{corollary}[theorem]{Corollary}
\newtheorem{lemma}[theorem]{Lemma}

\numberwithin{equation}{section}

\newcommand{\pf}{\textbf{Proof}\, }

\DeclareMathOperator{\Op}{Op}%
\DeclareMathOperator{\Smb}{Smb}%
\DeclareMathOperator{\supp}{supp}%
\DeclareMathOperator{\Imm}{Im}%

\renewcommand{\r}{\mathbb{R}}%
\newcommand{\cs}{\mathcal{S}}
\newcommand{\func}{\text}%
\newcommand{\tq}{\, \big| \, }%

\newcommand{\interior} {\selectlanguage{greek} \text{\bf i}} 
\newcommand{\ee} {\selectlanguage{greek} \text{\bf e}}

 \author{Marc Troyanov}
 \title{On the Hodge decomposition in  $\mathbb{R}^n$}
 
\begin{document}
 
\maketitle

\begin{abstract}
We prove a version of the $L^p$-Hodge decomposition for differential forms in Euclidean space and a generalization to the class of Lizorkin currents.  
Using these tools, we also compute the  $L_{q,p}-$cohomology of $\mathbb{R}^n$.
\end{abstract}

\maketitle
  
\section{Introduction}

The classical  \emph{Hodge decomposition Theorem} is a fundamental result in differential geometry. It says that on a compact Riemannian manifold without boundary, any $L^2$-differential form can be uniquely  decomposed as the sum of an exact form plus a coexact  and  a harmonic form. The standard proof is obtained by constructing an inverse $G$ to the Laplacian $\Delta$ on the $L^2$-orthogonal complement to the space of harmonic forms (see e.g. \cite[chap. 5.8]{taylor96}). The map $G$ is called the \emph{Green operator}. 
In 1995, Chad Scott proved that this decomposition also holds for $L^p$-differential  forms on a closed manifold \cite{scott95}.
For non compact, complete Riemannian manifolds, the $L^2$-Hodge decomposition holds under some technical hypothesis (which can be written in the language of $L^2$-cohomology, see \cite[Theorem 14.3]{GT2006}), but the situation for $L^p$-forms is still an open problem.

\medskip

The main goal of the present paper is to provide a rigorous statement and proof of the  Hodge decomposition theorem for $L^p$-differential forms in Euclidean space. Our proof uses standard techniques from Fourier analysis, but also the more specialized notion of \emph{Lizorkin currents} (see  sections 4 and 8.2  below and the books \cite{rubin96,samko02}). In order to keep the paper readable for non specialists
we give all the necessary background starting from the classic notion of tempered distribution, which are dual objects to rapidly decreasing smooth functions.

\medskip

Our techniques also provide a Hodge decomposition for \emph{temperate currents}, which are differential forms on $\mathbb{R}^n$ with coefficients in the space of temperate distributions  (besides a short mention in the last section of the book \cite{schwartz98}, temperate currents seem to have  been left aside in the literature).
Let us denote  by $\mathcal{S}'(\r^n, \Lambda^k)$  the topological vector
space of temperate currents of degree $k$, we then have the following theorem:

\begin{theorem} There is an   exact sequence
 \begin{equation}\label{ex.seq1}
 0\rightarrow \mathcal{H}(\r^n, \Lambda^k) \rightarrow \mathcal{S}'(\r^n,  \Lambda^k)\overset{\Delta}{\rightarrow }\mathcal{S}'(\r^n, \Lambda^k)\rightarrow 0,
\end{equation}%
where  $\mathcal{H}(\r^n, \Lambda^k)$ is the space of differential forms on $\r^n$ whose coefficients are  harmonic polynomials. 
\end{theorem}

This theorem is contained in Corollary \ref{cor.AA} below.
We also have an exact sequence 
 \begin{equation}\label{ex.seq1.1}
 0\rightarrow \mathcal{H}(\r^n, \Lambda^k) \rightarrow \mathcal{P}(\r^n,  \Lambda^k)\overset{\Delta}{\rightarrow }\mathcal{P}(\r^n, \Lambda^k)\rightarrow 0,
\end{equation}%
where $\mathcal{P}(\r^n, \Lambda^k)$ is the space of differential forms of degree $k$ on 
$\mathbb{R}^n$ with polynomial coefficients.

\medskip

It follows  from the sequence (\ref{ex.seq1}) and the identity $\Delta = (d\delta + \delta d)$,
that any $\theta\in\mathcal{S}'(\r^n, \Lambda^k)$ can be written as
\begin{equation}\label{hodge.dec1}
 \theta = \Delta \omega = d(\delta\omega)  + \delta(d\omega)
\end{equation}%
for some $\omega \in\mathcal{S}'(\r^n, \Lambda^k)$.
In particular,  any temperate current  $\theta\in\mathcal{S}'(\r^n, \Lambda^k)$ can be decomposed as the sum of an exact current plus a coexact current, this is the \emph{Hodge-Kodaira decomposition} for such currents. A similar statement holds for polynomial differential forms.

\medskip

If one formally introduces the operators 
$$ U = \delta \circ \Delta^{-1} , \quad U^* =  d\circ \Delta^{-1},$$
then the Hodge decomposition (\ref{hodge.dec1}) writes
\begin{equation}\label{hodge.dec2}
 \theta  = d(U\theta)  + \delta(U^*\theta).
\end{equation}%

Of course, due to the kernel  $\mathcal{H}(\r^n, \Lambda^k)$  in the exact sequence (\ref{ex.seq1}),
the operator $\Delta^{-1}$, as well as  $U$ and $ U^*$ are  not really well defined. 
However if one restricts our attention to  the space $L^p(\r^n,\Lambda^k)
\subset  \mathcal{S}' (\r^n,\Lambda^k)$ of 
differential forms with coefficients in $L^p(\r^n)$, then 
$$
 L^p(\r^n,\Lambda^k) \cap
 \mathcal{H}(\r^n,\Lambda^k) = \{ 0 \}
$$
for any $1\leq p < \infty$,  because $L^p(\r^n)$ does not contains any non zero polynomials.
The Laplacian $\Delta :  L^p(\r^n,\Lambda^k) \to \mathcal{S}' (\r^n,\Lambda^k)$ is thus  injective
and the operators $U$ and $ U^*$ can be properly defined on appropriate 
subspaces of $ \mathcal{S}' (\r^n,\Lambda^k)$.

\bigskip

We can now state the Hodge-Kodaira decomposition for the space
$L^p(\r^n,\Lambda^k)$:

\begin{theorem}\label{th.hodgeLp} Let $1<p<\infty$.
The space $L^p(\r^n,\Lambda^k)$ admits the following direct sum decomposition
\begin{equation}\label{eq.HKDLP}
 L^p(\r^n,\Lambda^k) = EL^p(\r^n,\Lambda^k) \oplus E^*L^p(\r^n,\Lambda^k),
\end{equation}
where  $EL^p(\r^n,\Lambda^k) = L^p(\r^n,\Lambda^k)\cap d\mathcal{S}'(\r^n,\Lambda^{k-1} T^*\r^n)$ is the space of exact currents belonging to $L^p$ and $E^*L^p(\r^n,\Lambda^k) = L^p(\r^n,\Lambda^k)\cap \delta\mathcal{S}'(\r^n,\Lambda^{k+1} T^*\r^n)$ is the space of coexact currents belonging to $L^p$.
Furthermore:
\begin{enumerate}[i.)]
  \item $EL^p(\r^n,\Lambda^k)$ and $E^*L^p(\r^n,\Lambda^k)$ are
  closed subspaces;
  \item the projections $E : L^p(\r^n,\Lambda^k)\to EL^p(\r^n,\Lambda^k)$
 and  $E^* : L^p(\r^n,\Lambda^k)\to E^*L^p(\r^n,\Lambda^k)$ are
 bounded operators;
 \item  these operators  satisfy
 $$
  E^2 = E, \qquad {E^*}^2 = E^*, \qquad  E + E^* = \mathrm{Id};
 $$
  \item the projection $E$ is self-adjoint, meaning that if $\frac{1}{p}+\frac{1}{q}=1$,
  then 
    $$\langle E\theta , \phi\rangle = \langle \theta , E\phi\rangle,$$
 for any   $ \theta \in L^p(\r^n,\Lambda^k)$ and $\phi \in L^q(\r^n,\Lambda^k)$ (where $\langle\theta , \phi\rangle = \int_{\mathbb{R}^n}  \theta \wedge \star \phi$);
 \item the same property holds for $E^*$.
\end{enumerate}
\end{theorem}

\medskip

We proove this result in section \ref{sec.LpTheory}.

\bigskip

\begin{corollary}
If $\frac{1}{p}+\frac{1}{q}=1$, then $E^*L^q(\r^n,\Lambda^k)$ and $EL^p(\r^n,\Lambda^k)$ are orthogonal, meaning that 
$$
 \langle\theta , \phi\rangle = \int_{\mathbb{R}^n}  \theta \wedge \star     \phi =  0,
$$
for any   $ \theta \in EL^p(\r^n,\Lambda^k)$ and $\phi \in E^*L^q(\r^n,\Lambda^k)$.
\end{corollary}

\textbf{Proof.}  Since $E$ and $E^*$ are projectors on complementary subspaces, we have $E\theta = \theta$ for any $ \theta \in EL^p(\r^n,\Lambda^k)$ and $E\phi = 0$ for any $\phi \in E^*L^q(\r^n,\Lambda^k)$, hence $\langle \theta | \phi\rangle = \langle E\theta | \phi\rangle  =
\langle \theta | E\phi\rangle  = \langle \theta | 0\rangle  =  0.$

\qed

\bigskip

The previous Theorem implies that any differential form $\theta \in L^p(\r^n,\Lambda^k)$ admits a Hodge-Kodaira decomposition $\theta = d\alpha + \delta \beta$, 
where  $d\alpha= E\theta$ and $\delta \beta = E^*\theta$ belong to $L^p(\r^n,\Lambda^k)$. The forms $\alpha$ and $\beta$ are in general just temperate distributions, but more can be said if  $1<p<n$:

\begin{theorem}\label{th.hodgeLqp}  
Let $1<p<n$ and $q=\frac{np}{n-p}$. There are bounded linear operators
$$
U^*:L^{p}(\r^n,\Lambda ^{k})\rightarrow L^{q}(\r^n,\Lambda ^{k-1})
\quad \text{and} \quad
U:L^{p} (\r^n,\Lambda ^{k})\rightarrow L^{q}(\r^n,\Lambda ^{k+1}),
$$
such that 
$
E =d\circ U:L^{p}(\r^n,\Lambda ^{k})\rightarrow EL^{p}(\r^n,\Lambda ^{k})
$ and 
$E^* =\delta \circ U^* :L^{p}(\r^n,\Lambda ^{k})\rightarrow E^*L^{p}(\r^n,\Lambda ^{k})$. 
In particular,  any differential form $\theta
\in L^{p}(\r^n,\Lambda ^{k})$ can be uniquely decomposed as a sum of an exact
form $d\alpha $ plus a co-exact form $d\beta $ with $\alpha = U\theta \in
L^{q}(\r^n,\Lambda^{k-1})$ and $\beta = U^*\theta \in L^{q}(\r^n,\Lambda^{k+1})$:
\begin{equation} \label{dec.1}
 \theta = E\, \theta  + E^*\, \theta  = d(U\theta)  + \delta(U^*\theta) =d\alpha +\delta \beta .  
\end{equation}%
\end{theorem}

This theorem is also proved in Section \ref{sec.LpTheory}. We will also see that such a decomposition exists
only if $1<p<n$ and $q=\frac{np}{n-p}$.

\bigskip

As said before, our results are proved using various known facts from harmonic analysis and symbolic calculus in the
context of differential forms and currents. We present all necessary notions in a  self contained way and the paper is organized as follows: In section 2, we recall some basic facts about the space $\mathcal{S}'$ of tempered disitribution, the Fourier transform and the notion of convolution. In section 3 we recall the characterization of polynomials as elements  in 
$\mathcal{S}'$ annihilating some power of the Laplacian. In section 4 we introduce the Lizorkin distributions which are tempered distributions modulo the space of polynomials. In section 5, we develop some symbolic calculus for Lizorkin distributions and apply it to the Riesz potential and Riesz transform and in section 6 we recall some basic facts from $L^p$-harmonic analysis.

In sections 2--6, only functions are investigated. In section 7, we recall some basic facts about differential forms in $\r^n$ and we
prove all the results stated in the present introduction. Section 8 is devoted to some applications of the previous results,
in particular we prove three fundamental inequalities concerning differential forms in $\r^n$ and we give a necessary and sufficient 
condition for the vanishing of the $L_{qp}$-cohomology of $\r^n$. The paper ends with a technical appendix devoted to a
calculation of the Fourier transform of the Riesz kernel.

\section{The space of tempered distributions}

We will work with the Fourier transform of tempered distributions as they are developed e.g. in  \cite{schwartz98,strichartz94,taylor96}.  Recall that the  Schwartz space
$\mathcal{S}=\mathcal{S}(\mathbb{R}^{n})$ of rapidely decreasing
functions is the space of smooth functions $f:\mathbb{R}^{n}\rightarrow \mathbb{C}$ such that 
\begin{equation*}
  [f]_{m,\alpha} = \|(1+|x|)^m\partial ^{\alpha }f \|_{L^{\infty }(\mathbb{R}^{n})}
  < \infty
\end{equation*}
for all $m\in \mathbb{N}$ and  all multi-indices $\alpha\in
\mathbb{N}^n$. This is a Frechet space for the
topology induced by the collection of all semi-norms $[ \, \cdot \, 
]_{m,\alpha}$, it is  dense in $L^p(\r^n)$ for any $1\leq p
< \infty$ and it is also a pre-Hilbert space for the inner
product $(f|g)=\left\langle f,\overline{g}\right\rangle$ where
\begin{equation}\label{form.eval}
 \left\langle f,g\right\rangle =\int_{\r^n} f(x)g(x)  dx.
\end{equation}
Recall also that $\cs$ is an algebra for the multiplication and
for the convolution product, it is closed under translation,
differentiation and multiplication by polynomials. 
\medskip

Of basic importance is the fact  that the Fourier transform\footnote{There are different conventions for this definition, this affects some constants in the following formulas. Here, we follow \cite{strichartz94}.}
\[
\mathcal{F}(f)(\xi )=\widehat{f}(\xi
)=\int_{\mathbb{R}^{n}}f(x)e^{ix\cdot \xi }dx
\]%
is an isomorphism $\mathcal{F}:\mathcal{S}\rightarrow
\mathcal{S}$, with inverse
\[
 \mathcal{F}^{-1}(g)(x)=\check{g}(x) =
 \frac{1}{(2\pi)^{n}}\int_{\mathbb{R}^{n}}g(\xi)e^{-ix\cdot\xi }d\xi.
\]
Some of the basic properties of the Fourier transform
are\begin{enumerate}[i)]
\item $\mathcal{F}(f\ast g)=\mathcal{F}(f)\cdot \mathcal{F}(g)$ \ (here and below $\ast$ is the convolution product);

\item $\ \mathcal{F}(f\cdot g)=\frac{1}{(2\pi )^{n}}\mathcal{F}(f)\ast \mathcal{F}(g);$

\item $\mathcal{F}(\partial _{j}f)(\xi )=-i\xi _{j}\cdot \mathcal{F}(f)(\xi );$

\item $\mathcal{F}(\bar{g})=(2\pi )^{n}\overline{\mathcal{F}^{-1}(g)}$;

\item $\mathcal{F}\left( f\circ A\right) =\frac{1}{|\det
A|}\mathcal{F}(f)\circ \left( A^{-1}\right) ^{t}$ for any $A\in
GL_{n}$\bigskip $(\mathbb{R)}$.
\end{enumerate}
From Fubini's Theorem, we have
\begin{equation}\label{dualF}
\left\langle \mathcal{F}f,g\right\rangle =\left\langle
f,\mathcal{F}g\right\rangle 
= \int_{\mathbb{R}^{n}}\int_{\mathbb{R}^{n}}f(x)g(y)e^{ixy}dxdy.
\end{equation}
This identity can also be written as
\begin{equation}
\left(\mathcal{F}f | g\right) =
(2\pi )^{n}\left(f | \mathcal{F}^{-1}g\right)
 \quad   \textrm{or}   \quad
\left(\mathcal{F}f|\mathcal{F}h\right) =
(2\pi )^{n}\left(f|h\right)
\end{equation}
(just set $h=\mathcal{F}^{-1}g$ in the previous identity). The
latter formula is  the   \emph{Parseval-Plancherel identity}.

\medskip 

The topological dual of $\mathcal{S}$  is called the space of \emph{tempered distributions} and
is denoted by $\mathcal{S}'$, and if $w\in \cs'$ and $f\in \cs$, the
evaluation of $w$ on $f$ will be denoted by
\begin{equation*}
  \langle w,f\rangle \in \mathbb{C}.
\end{equation*}
Any measurable function $f$ such that $|f(x)| \leq C (1+|x|^m)$
for some $m>0$ and any function in $f\in L^p(\r^n)$ defines a
tempered distributions\footnote{More generally, a complex Borel
measure  $\mu$ on $\mathbb{R}^{n}$  belongs to $\cs'$ if and only
if $|\mu (B(0,R))|\leq C\cdot (1+R)^{N}$ for some $N\in
\mathbb{Z}$ and all $R>0$.} by the formula (\ref{form.eval}).
Distributions with compact support also  belong to $\cs'$.
 
\medskip 

The space  $\mathcal{S}'$ is a complete locally convex 
topological vector space when equipped with the weak*
topology, i.e. the smallest topology for which  the linear form
\begin{equation*}
 w \mapsto  \langle w,\varphi\rangle 
 \end{equation*}
is continuous for any $\varphi \in \mathcal{S}$ (note that  $\mathcal{S}'$ is not a Frechet space).

\begin{lemma}\label{w*closed}
 If $A\subset \r^n$ is a non empty closed subset, then 
 $$
    \mathcal{S}'_A = \{ w \in  \mathcal{S}' \tq  \supp (w) \subset A \}
 $$
 is a closed subset in $\mathcal{S}'$.
\end{lemma}

\pf  Suppose that $w_0 \notin   \mathcal{S}'_A$, then, by definition, there
exists a function $\varphi \in   \mathcal{S}$ such that $\supp (\varphi) \cap A = \emptyset$ and $s =  \langle w_0,\varphi \rangle > 0$. \
Consider now the set $\mathcal{W} \subset \mathcal{S}'$ defined by
$$
\mathcal{W} = \left\{ w \in  \mathcal{S}' \tq  \langle w,\varphi \rangle > \frac{s}{2} \right\}.
$$
By definition of the weak* topology, $\mathcal{W}$ is open in  $\mathcal{S}'$.
It is clear that $\mathcal{W}\cap \mathcal{S}'_A=\emptyset$. We have thus 
found, for any $w_0 \notin   \mathcal{S}'_A$, an open set such that 
$$w_0 \in \mathcal{W} \subset  \mathcal{S}' \setminus \mathcal{S}'_A.$$
This  means that the complement of $\mathcal{S}'_A$ is an open
subset in $\mathcal{S}'$.

\qed

\bigskip

The differential operator $\partial_i$ acts continuously on $\cs'$ by duality:
\begin{equation*}
  \langle\partial_i w,f\rangle = - \langle w,\partial_i f\rangle.
\end{equation*}

We can also define the Fourier transform by
\begin{equation*}
  \langle\mathcal{F} w,f\rangle = \langle w,\mathcal{F}f\rangle
\end{equation*}
and its inverse by
\begin{equation*}
  \langle\mathcal{F}^{-1} w,f\rangle = \langle w,\mathcal{F}^{-1} f\rangle.
\end{equation*}
These are continuous isomorhisms $\mathcal{F} , \mathcal{F}^{-1}
:\mathcal{S}'\rightarrow \mathcal{S}'$ which are inverse to each
other.
Some important examples of Fourier transforms are  
$$
 \mathcal{F}(\mathrm{e}^{-\frac{|x|^2}{2}})(\xi) = (2\pi)^{\frac{n}{2}}\mathrm{e}^{-\frac{|\xi|^2}{2}}, \qquad
  \mathcal{F}(1) = (2\pi)^n\delta_0, \qquad
  \mathcal{F}(\delta_0) = 1,
$$
where $\delta_0 \in \cs'$ is the  Dirac measure.

\medskip

The convolution of two tempered distributions is in general not
defined, but we can define a convolution product
\begin{equation*}
  * : \cs \times \cs' \to \cs'
\end{equation*}
by the formula
\begin{equation*}
  \langle f* w,g\rangle = \langle w,\tilde{g}* f\rangle , 
  \end{equation*}
where $w \in \cs'$ and $f,g\in \cs$. Here 
 $\tilde{f}(x) = f(-x)$. Observe that this formula is
consistent with  Fubini theorem in the case $w\in \cs$. The Dirac
measure $\delta_0 \in \cs'$ is the convolution identity in the sense that
$$ f* \delta_0 = f.$$
for all $f\in \cs$.

\section{The Laplacian and Polynomials}

Let us denote by $\mathcal{P}$ the space of all polynomials
$P:\mathbb{R}^{n}\rightarrow \mathbb{C}$. It is a subspace of  $\mathcal{S}$  and it has the following important characterization (see  \cite[Proposition 4.5]{taylor96}):

\begin{proposition}\label{prop.FP}
A tempered distribution $f\in \mathcal{S}'$ is a polynomial if and
only if the support of its Fourier transform is contained in
$\{0\}$:
\begin{equation*}
\mathcal{P}=\left\{ f\in \mathcal{S}' \tq \supp\hat{f}\subset
\{0\}\right\} .
\end{equation*}
\end{proposition}

\begin{corollary}
 $\mathcal{P}$ is a closed subspace of  $\mathcal{S}'$.
\end{corollary}

\pf This follows directly from the previous Proposition and 
Lemma \ref{w*closed}.

\qed

\medskip

The Laplacian on $\mathbb{R}^{n}$\ is the partial differential
operator $\Delta = - \sum_{j=1}^{n}\partial _{j}^{2}$. For a
distribution $w\in \mathcal{S}'$, we define $\Delta w\in \mathcal{S}'$
by 
\begin{equation*}
\left\langle \Delta w,\varphi \right\rangle =\left\langle
w,\Delta \varphi \right\rangle
\end{equation*}%
for any $\varphi \in \mathcal{S}$. 
The relation with the Fourier transform is given by
\begin{equation*}
\mathcal{F}(\Delta w)(\xi ) = |\xi |^{2}\mathcal{F}(w).
\end{equation*}

A distribution $w\in \mathcal{S}' $, is called \emph{harmonic} if $\Delta w=0$ and
we denote by $\mathcal{H}$ the space of harmonic tempered distributions,
i.e. the kernel of $\Delta$ :
$$
 \mathcal{H} = \{w\in \mathcal{S}' \tq  \Delta w=0 \}.
$$

\bigskip

\begin{proposition} \label{prop.kerdelm}
A tempered distribution $f\in \mathcal{S}'$ is a polynomial if and
only if $\Delta^{m}f=0$ for some $m\in \mathbb{N}$.
\end{proposition}

\bigskip

\textbf{Proof} It is obvious that if $f\in \mathcal{P}$ is a
polynomial of degree $m$, then $\Delta ^{m}f=0$. Conversely, if
$\Delta ^{m}f=0$, then
\begin{equation*}
0=\mathcal{F}(\Delta ^{m}f)(\xi )=(-1)^{m}|\xi |^{2m}\hat{f}(\xi),
\end{equation*}%
hence  $\supp\hat{f}\subset \{0\}$.

\qed

\bigskip

We just proved that $\mathcal{H} = \ker \Delta \subset \mathcal{P}$. A consequence of this 
result is the following generalization of Liouville's theorem:

\begin{corollary}
\begin{equation*}
L^{\infty }(\mathbb{R}^{n})\cap \ker \Delta =\mathbb{R\qquad
}\mbox{and}\qquad L^{p}(\mathbb{R}^{n})\cap \ker \Delta
=\{0\}\mathbb{\quad }\text{for }1\leq p<\infty .
\end{equation*}%
\end{corollary}

\qed

Observe however that not every globally defined harmonic function in
$\mathbb{R}^{n}$ is a polynomial, for instance the function
$h(x)=\sin (x_{1})\sinh (x_{2})$ is harmonic. Of course $h\notin
\mathcal{S}'$.

\bigskip

\begin{theorem}\label{Th.InvPoly}
The Laplacian $\Delta :\mathcal{P}\rightarrow \mathcal{P}$\ is
surjective;\ we thus have an exact sequence
\begin{equation*}
0\rightarrow \mathcal{H}\rightarrow \mathcal{P}\overset{\Delta
}{\rightarrow }\mathcal{P}\rightarrow 0.
\end{equation*}%
\end{theorem}

\bigskip

\textbf{Proof} A computation show that if $m\in \mathbb{N}$ and
$h\in \mathcal{P}$ is a homogenous function of degree
$\nu $ (i.e. $h(tx)=t^{\nu}h(x)$ for $t>0$), then
\begin{equation}\label{Deltapoly}
\Delta (|x|^{2m+2}h(x))=|x|^{2m+2}\Delta h(x)+c_{n,m,\nu }|x|^{2m}h(x)
\end{equation}%
where $c_{n,m,\nu } = 2(m+1)(2m+2\nu +n)$ (use Euler's Formula for homogenous functions:   $\sum  x_i\partial_i h = \nu \cdot h$).
On the other hand, a basic result about polynomials (see \cite{AR}) says  that any $f\in \mathcal{P}$ can be
written as a finite sum
\begin{equation}\label{frmule.h}
 f(x)=\sum_{m,\nu }a_{m,\nu }|x|^{2m}h_{m,\nu }(x)
\end{equation}%
where $h_{m,\nu }\in \mathcal{H}$ is a  homogenous polynomial of degree $\nu$
Now it is clear from (\ref{Deltapoly}) that \ $f=\Delta g$ with
\begin{equation}\label{frmule.g}
 g(x)=\sum_{m,\nu }\frac{a_{m,\nu
 }}{c_{n,m,\nu}}|x|^{2m+2}h_{m,\nu }(x).
\end{equation}%
The surjectivity of $\Delta :\mathcal{P}\rightarrow \mathcal{P}$ follows.

\qed

\bigskip

\textbf{Remark.}  The paper \cite{AR} gives an explicit procedure to compute the decomposition (\ref{frmule.h}),
the proof thus shows that the inverse Laplacian $\Delta^{-1} : \mathcal{P} \to \mathcal{P}$ given by  (\ref{frmule.g})
is algorithmically computable.

\section{The Lizorkin space and its Fourier image}

\textbf{Definition}  We introduce two subspaces  $\Phi$ and  $\Psi $ 
of $\mathcal{S}$ defined as follow:
$$
 \Phi = \bigcap_{m=0}^{\infty }\Delta ^{m}(\mathcal{S})
\quad  \text{ and }   \quad
\Psi = \left\{ \psi \in
\mathcal{S}:\mathbb{\partial }^{\mu }\psi (0)=0,\text{ for any } \mu
\in \mathbb{N}^{n}\right\} .
$$
The space $\Phi$ is called the \emph{Lizorkin space}, basic references on this space are
\cite{rubin96,samko02}. We shall see below that $\Psi$ is the Fourier dual of $\Phi$, that is  the image of $\Phi$ under the Fourier transform.

\begin{theorem}
The restriction of the Laplacian to the Lizorkin space is a
bijection $\Delta :\Phi \rightarrow \Phi $.
\end{theorem}

\textbf{Proof} The Laplacian is injective on $\mathcal{S}$
because $\ker \Delta \cap \mathcal{S}\subset \mathcal{P}\cap
\mathcal{S}=\{0\}$ by proposition \ref{prop.kerdelm}. 
To prove the surjectivity, consider an
arbitrary element $\varphi \in \Phi $. By definition, for any
$m\in \mathbb{N}$ there exists $g_{m}\in \mathcal{S}$ such that
$\Delta ^{m}g_{m}=\varphi $. Observe that $\Delta (\Delta
^{m}g_{m+1}-g_{1})=\varphi -\varphi =0$. Since $\Delta $ is
injective on $\mathcal{S}$, we have $g_{1}=\Delta ^{m}g_{m+1}\in
\Delta ^{m}(\mathcal{S})$. It follows that $g_{1}\in \Phi $ and
therefore $\varphi =\Delta g_1\in \Delta \Phi $.

\qed

\begin{proposition}\label{prop.charpsi}
For any rapidly decreasing function $\psi \in \mathcal{S}$,
the following conditions are equivalent:
\begin{enumerate}[(a)]
\item $\psi\in \Psi $;
\item $\mathbb{\partial }^{\mu }\psi(\xi)=o(|\xi |^{t})$ as $|\xi
|\rightarrow 0$ for any multi-indice $\mu \in \mathbb{N}^{n}$ and
any $t>0$;
\item $|\xi |^{-2m}\psi \in \mathcal{S}$ for any $m\in \mathbb{N}$.
\end{enumerate}
\end{proposition}

\textbf{Proof}  The implication (b)$\Rightarrow $(a) is obvious and (a)$\Rightarrow $(b) is clear by Taylor expansion.

\smallskip

To prove that (b)$\Rightarrow $(c), observe that condition (b),  together with the Leibniz rule, implies that the function $|\xi |^{-2m}\psi$  vanishes at the origin and is continuous as well as all its derivatives. It is then clear that $|\xi |^{-2m}\psi \in \mathcal{S}$. 

To prove  (c)$\Rightarrow $(a), observe that condition (c) says that
$\psi =|\xi |^{2m}\rho $ for some function $\rho \in \mathcal{S}$. 
By the Leibniz rule, we then have $\mathbb{\partial}^{m}\psi(0)=\partial ^{m}\left( |\xi |^{2m}\rho \right)
(0)=0$.

\qed

\begin{proposition}
We have $\mathcal{F}(\Phi) = \Psi$.
\end{proposition}

\textbf{Proof} \  For any\ $\varphi \in \Phi
$ and $m\in \mathbb{N}$ there exists $\varphi _{m}\in
\mathcal{S}$ such that $\Delta ^{m}\varphi _{m}=\varphi $. The
Fourier transform of this relation writes $\widehat{\varphi}=(-1)^{m}|\xi |^{2m}\widehat{\varphi _{m}}$, thus $|\xi
|^{-2m}\widehat{\varphi }=(-1)^{m}\widehat{\varphi _{m}}\in
\mathcal{S}$ for any integer $m$ and it follows from condition (c) in the previous
 proposition that  $\widehat{\varphi } \in \Psi$, hence $\mathcal{F}(\Phi) 
 \subset \Psi$.

 \smallskip

To prove the opposite inclusion, we consider a function $\psi \in \Psi$. 
Using again condition (c) in the previous proposition, we know that for any $m\in \mathbb{N}$, 
we can write $\psi =|\xi |^{2m}\psi_m$ for some
function $\psi_m \in \mathcal{S}$. We then have 
$$
 \mathcal{F}^{-1}(\psi) = \mathcal{F}^{-1}(|\xi |^{2m}\psi_m)
 = (-1)^{m}\Delta ^{m}\left( \mathcal{F}^{-1}(\psi_m) \right),
$$
hence $\mathcal{F}^{-1}(\psi)  \in 
 \bigcap_{m=0}^{\infty }\Delta ^{m}(\mathcal{S}) = \Phi$.
 
 \qed

\begin{corollary}
For any   $\varphi \in \mathcal{S}$, we have 
$$
\varphi \in \Phi  \ \Leftrightarrow  \
\left\langle P,\varphi \right\rangle =0 \ \text{ for any polynomial } P\in
\mathcal{P}.
$$
\end{corollary}

\textbf{Proof} \  
For any $\varphi \in \mathcal{S}$ and
$\mu \in \mathbb{N}^{n}$, we have
\begin{equation*}
\left\langle x^{\mu },\varphi \right\rangle
=\int_{\mathbb{R}^{n}}x^{\mu }\varphi (x)dx=i^{-|\mu
|}\int_{\mathbb{R}^{n}}(ix)^{\mu }\varphi (x)e^{-ix\cdot
0}dx=i^{-|\mu |}{\partial }^{\mu }\hat{\varphi}(0).
\end{equation*}%
Thus $\left\langle P,\varphi \right\rangle =0$  for any polynomial
if and only if ${\partial }^{\mu }\hat{\varphi}(0) = 0$ for any $\mu \in \mathbb{N}^{n}$,
i.e. if $\hat{\varphi} \in \Psi$ and we conclude by the previous Proposition.

\qed

\bigskip
\begin{proposition}
$\Psi $ is a closed ideal of $\mathcal{S}$.
\end{proposition}

\textbf{Proof} It is clear from the Leibniz rule that if $\psi
\in \Psi $ and $f\in \mathcal{S}$, then $\mathbb{\partial }^{\mu}
(f\psi )(0)=0$ for any $\mu \in \mathbb{N}^{n}$, hence $\Psi
\subset \mathcal{S}$ is an ideal. To show that $\Psi \subset
\mathcal{S}$ is closed, let us consider a sequence $\{\psi
_{j}\}\subset \Psi $ converging to $\psi \in \mathcal{S}$. This
means that for any $\mu \in \mathbb{N}^{n}$, $m\in \mathbb{N}$,
we have $\sup_{x\in \mathbb{R}^{n}}(1+|x|)^{m}\left|
\mathbb{\partial }^{\mu }\psi - \mathbb{\partial }^{\mu }\psi _{j}\right| \rightarrow 0$ as
$j\rightarrow \infty $. But then $\ds \mathbb{\partial }^{\mu }\psi
(0)=\lim_{j\rightarrow \infty }\mathbb{\partial }^{\mu }\psi
_{j}(0)=0$.

\qed

\bigskip

Since $\mathcal{F}:\mathcal{S}\rightarrow \mathcal{S}$ is a linear
homeomorphism sending $\Phi $ to $\Psi $, we immediately conclude
that
\begin{corollary}
$\Phi $ is a closed subspace of $\mathcal{S}$ and it is an ideal
for the convolution.
\end{corollary}

\medskip

The statement that $\Phi \subset \mathcal{S}$ is a convolution
ideal means that if  $\varphi \in \Phi $ and $f\in \mathcal{S}$,
then $\varphi \ast f=f\ast \varphi \in \Phi $. This can also be
seen directly from the definition of $\Phi $, indeed, if $\varphi
\in \Phi $, then for any $m\in \mathbb{N}$ there exists $g_{m}\in
\mathcal{S}$ such that $\Delta ^{m}g_{m}=\varphi $ and we have
\begin{equation*}
 \Delta ^{m}(f\ast g_{m})=f\ast (\Delta ^{m}g_{m})=f\ast \varphi .
\end{equation*}%
Thus $f\ast \varphi \in \Delta ^{m}(\mathcal{S})$, for any $m$,
i.e. $f\ast \varphi \in \Phi $.

\begin{corollary}
$\Delta :\Phi \rightarrow \Phi $ is a homeomorphism.
\end{corollary}

\textbf{Proof} We already know that $\Delta :\Phi \rightarrow
\Phi $ is bijective. The inverse $\Delta ^{-1}:\Phi \rightarrow
\Phi $ is given by the formula
\begin{equation*}\label{inv.lap}
 \Delta ^{-1}(\varphi) = \mathcal{F}^{-1}\left(|\xi|^{-2}\mathcal{F}\varphi \right).
\end{equation*}
Since the map $\psi \to |\xi|^{2}\psi$ is clearly a self-homeomorphism of $\Psi$,
we obtain the continuity of $\Delta ^{-1}:\Phi \rightarrow
\Phi $.

\qed

\bigskip

\begin{proposition}
(A) The topological dual $\Phi '$ of $\Phi $ is the quotient of
the space of tempered distribution modulo the polynomials
\begin{equation*}
 \Phi '=\mathcal{S}'/\mathcal{P}.
\end{equation*}
(B) The topological dual $\Psi '$ of $\Psi $ is the quotient of
the space of tempered distribution modulo the Fourier transforms
of polynomials
\begin{equation*}
 \Psi '=\mathcal{S}'/\mathcal{F(P)}.
\end{equation*}
\end{proposition}

\pf  The closed subspace 
$\mathcal{P}\subset \mathcal{S}'$ coincides with $\Phi^{\perp} = \{ w \in \mathcal{S}'
\tq w(\Phi) = 0\}$ and $\Psi^{\perp} = \mathcal{F} (\mathcal{P})$ . 
The Proposition follows now from standard results from functional analysis (see e.g.  \cite[chap. V, th. 2.3]{conway}).

\qed

\bigskip

An element $w\in \Phi '$ is thus represented by a tempered
distribution which is only well defined up to a polynomial.
 The Fourier transform $\mathcal{F}:\mathcal{S}'\rightarrow
\mathcal{S}'$ gives an isomorphism between these quotients which
we continue to denoted by $\mathcal{F}:\Phi '\rightarrow \Psi '$. We have
\begin{equation*}
\langle\mathcal{F}w,\varphi\rangle = \langle
w,\mathcal{F}\varphi\rangle
\end{equation*}
for any $w\in \Phi '$ and $\varphi \in \Phi $.

\section{Some symbolic calculus}

\subsection{Operators on $\Psi'$  and multipliers}

In this section, we study the operators $M:\Psi '\rightarrow
\Psi'$, which can be represented by a multiplication.

\medskip

\textbf{Definitions \,  1)} By an \emph{operator }
$M:\Psi'\rightarrow \Psi'$, we mean a continuous linear map.
Concretely, an operator associates to an element $w\in S'$ another
tempered distribution $Mw \in S'$ which is well defined modulo
$\mathcal{F(P)}$.
The linearity means that
$M(a_{1}w_{1}+a_{2}w_{2})=a_{1}M(w_{1})+a_{2}M(w_{2})$ modulo
$\mathcal{F(P)}$ for any $a_{1},a_{2}\in \mathbb{C}$,
$w_{1},w_{2}\in \mathcal{S}'$ and the  continuity means that $\left\langle w_{i},\psi \right\rangle
\rightarrow \left\langle w,\psi \right\rangle $ for any $\psi \in
\Psi $ implies $\left\langle Mw_{i},\varphi \right\rangle
\rightarrow \left\langle Mw,\varphi \right\rangle $. If $M$ has a
continuous inverse, then we say that it is an \emph{isomorphism}.

\medskip

\textbf{2)} We denote by $\Op(\Psi ')$ the algebra of all
operators $\Psi '\rightarrow \Psi '$.

\bigskip

We will discuss a special class of operators on $\Psi$, obtained
by multiplication with a suitable function, which we now
introduce:

\medskip

\textbf{Definition} Let $\mathcal{M}_{\Psi '}$ be the space of all
functions  $\sigma \in C^{\infty }(\mathbb{R}^{n}\setminus
\{0\},\mathbb{C)}$ such that for any multi-index $\mu \in
\mathbb{N}^{n}$, there exists constants $m\in \mathbb{N}$ and
$C>0$ with
\begin{equation*}
  \left| \partial ^{\mu }\sigma (\xi )\right|
  \leq C\left( |\xi|^{m}+ |\xi |^{-m}\right)
\end{equation*}
for any $\xi \in \mathbb{R}^n\setminus \{0\}$. An element of $\mathcal{M}_{\Psi '}$
is called a ${\Psi '}$-\emph{multiplier}.

\medskip

It is that clear that $\mathcal{S}\subset \mathcal{M}_{\Psi'}$
and $\mathcal{P}\subset \mathcal{M}_{\Psi '}$, other typical
elements of $\mathcal{M}_{\Psi '}$ are the
functions $\log |\xi|$ and $|\xi |^{\alpha }$ (for any $\alpha \in
\mathbb{C}$). Observe also that $\mathcal{M}_{\Psi '}$ is a
commutative algebra.

\medskip

The units in $\mathcal{M}_{\Psi '}$, i.e. the group of invertible
elements, will be denoted by $\mathcal{UM}_{\Psi '}$, hence
$$
\mathcal{UM}_{\Psi '} = \{\sigma \in \mathcal{M}_{\Psi '} \tq
\frac{1}{\sigma} \in \mathcal{M}_{\Psi '}\}.
$$

\medskip

Elements in $\mathcal{M}_{\Psi '}$ are not tempered
distributions, however, we have the following important lemma:

\begin{lemma}\label{psistable} 
$\Psi$ is a module over the algebra
$\mathcal{M}_{\Psi '}$, that is for any $\sigma \in \mathcal{M}_{\Psi '}$ and any $\psi\in \Psi$,
we have $\sigma \cdot \psi \in \Psi$.
\end{lemma}

\textbf{Proof}  By Proposition \ref{prop.charpsi}, we
know that an element $\psi \in \mathcal{S}$ belongs to $\Psi$ if
and only if $\mathbb{\partial }^{\mu }\psi(\xi)=o(|\xi |^{t})$ as
$|\xi |\rightarrow 0$ for any multi-indice $\mu \in
\mathbb{N}^{n}$ and any $t>0$. The proof of the lemma follows now
easily from the Leibniz rule.

\qed

By duality, we can now associate to any $\sigma \in
\mathcal{M}_{\Psi '}$ an operator $M_{\sigma } \in \Op (\Psi')$
defined by
\begin{equation}
    \langle M_{\sigma } g , \psi \rangle = \langle g , \sigma\psi \rangle
\end{equation}
for any $g\in \Psi'$ and  $\psi \in \Psi$.

\begin{lemma}
This correspondence defines a map
\begin{eqnarray*}
M  : \mathcal{M}_{\Psi '} &\rightarrow & \Op (\Psi') \\
     \sigma   &\mapsto &  M_{\sigma }
\end{eqnarray*}%
which is a continuous homomorphism of algebras. 
In particular $M_{\sigma_1 \sigma_2} = M_{\sigma_1}\circ
M_{\sigma_2}$ and $M_{\sigma }$ is invertible if and only if
$\sigma \in \mathcal{UM}_{\Psi '}$
\end{lemma}
The proof is elementary. 

\qed

\bigskip

\textbf{Definition} An operator $M_{\sigma } \in \Op (\Psi')$ of
this type is called a  \emph{multiplier in} $\Psi'$; the set of
those multipliers is denoted by $\mathcal{M}\Op (\Psi')$, it is a
commutative subalgebra of $\Op (\Psi')$.

\bigskip

Observe  that, by Lemma \ref{psistable},  $\Psi \subset \Psi'$ is
invariant under any multiplier in $\Psi'$ (i.e. $M_{\sigma}(\Psi)
\subset \Psi$ for any  $M_{\sigma } \in \mathcal{M}\Op (\Psi')$).
The converse is in fact also true.
The multiplication $M_{\sigma }(\psi )=\sigma \cdot \psi $ is a
continuous operator on $\Psi '$ if and only if $\sigma \in
\mathcal{M}_{\Psi '}$ (see \cite{samko77}).

\subsection{Operators on $\Phi'$ and their symbols}

In this section, we study the operators $T:\Phi '\rightarrow
\Phi'$, which can be represented on the Fourier side by a
multiplication. We already know the Laplacian:

\begin{equation*}
 \mathcal{F}(\Delta w) = |\xi |^{2}\mathcal{F}(w).
\end{equation*}%

\textbf{Definitions. \,  1)} An \emph{operator }$T:\Phi
'\rightarrow \Phi '$, is a continuous linear map, it associates to an element $w\in S'$ another tempered
distribution $Tw\in S'$ which is well defined up to a polynomial.

\medskip

\textbf{2)} We denote by $\Op(\Phi ')$ the algebra of all operators $\Phi '\rightarrow \Phi '$.

\medskip

An obvious consequence of the previous section is the ing

\begin{proposition}
For any $\sigma \in \mathcal{M}_{\Psi '}$, the map $T_{\sigma}  :
\Phi '\rightarrow \Phi '$ defined by
$$
T_{\sigma}=\mathcal{F}^{-1}\circ M_{\sigma }\circ \mathcal{F}
$$
belongs to $\Op(\Phi ')$. If $\sigma\in\mathcal{UM}_{\Psi '}$,
then $T_{\sigma }$ is an isomorphism of $\Phi '$.
\end{proposition}

\qed

\bigskip

\textbf{Definition.} Operators of this type are called
\emph{Fourier multipliers in} $\Phi'$. We denote the set of
those operators by $\mathcal{FM}\Op(\Phi ')$.

\medskip

If $T=T_{\sigma }= \mathcal{F}^{-1}\circ M_{\sigma }\circ
\mathcal{F} \in \mathcal{FM}\Op(\Phi ')$, then the function
$\sigma \in \mathcal{M}_{\Psi '}$ is the \emph{symbol} of
$T_{\sigma }$  and we write
$$\sigma = \Smb(T).$$

\bigskip

Thus, to say that $T \in \mathcal{FM}\Op(\Phi ')$ means that for
any Lizorkin distribution  $f\in \Phi'$ and any $\varphi \in
\Phi$, we have
$$
\langle Tf , \varphi \rangle = \langle f , \mathcal{F}\left(\sigma
\cdot \mathcal{F}^{-1}(\varphi)\right) \rangle,
$$
where $\sigma = \Smb(T)$.
Observe that $\mathcal{FM}Op(\Phi ')$ is a commutative algebra and
the  map $\mathcal{M}_{\Psi '}\rightarrow \mathcal{FM}\Op(\Phi ')$
given by $\sigma \mapsto T_{\sigma}$  is an isomorphism whose inverse is give by the symbol map:
\begin{eqnarray*}
\Smb: \mathcal{FM}\Op(\Phi ') &\rightarrow&
\mathcal{FM}_{\Psi '}.
\\
\sigma &\rightarrow& T_{\sigma}
\end{eqnarray*}

\medskip

\textbf{Examples : }
\begin{enumerate}[i)]
\item The symbol of the identity is $1$;
\item  $\Smb(T\circ U)= \Smb(T)\cdot \Smb(U)$;
\item The derivative $\partial _{j} \in \mathcal{FM}Op(\Phi ')$ and
$\Smb(\partial _{j})=-i\xi _{j}$;
\item  The symbol of the Laplacian is $\Smb(\Delta )= |\xi |^{2}$;
\item  More generally, $T$ is a partial differential operator with constant coefficients if
and only if $P=\Smb(T)\in \mathcal{P}$;
\item  If $T(w)=\varphi \ast w$ for some $\varphi \in \mathcal{S}$,
then $\Smb(T)=\widehat{\varphi }$.
\end{enumerate}

\medskip

Any operator $T_{\sigma }\in \mathcal{FM}\Op(\Phi ')$ is
self-adjoint in the following sense:

\medskip

\begin{proposition}
For any $T_{\sigma }\in \mathcal{FM}\Op(\Phi ')$ , we have
$T_{\sigma }(\Phi )\subset \Phi$ and
\begin{equation*}
 \left\langle T_{\sigma }w,\varphi \right\rangle =
 \left\langle w,T_{\sigma }\varphi \right\rangle
\end{equation*}
for all $w\in \Phi '$ and $\varphi \in \Phi $.
\end{proposition}
\textbf{Proof} The fact that $T_{\sigma }(\Phi )\subset \Phi$
follows from Lemma \ref{psistable} and we  have
\begin{eqnarray*}
\left\langle T_{\sigma }w,\varphi \right\rangle  &=&
\left\langle \mathcal{F}^{-1}\circ M_{\sigma }\circ \mathcal{F}(w),\varphi \right\rangle  \\
&=&(2\pi )^{-n}\left\langle M_{\sigma }\circ \mathcal{F}(w),\mathcal{F}(\varphi )\right\rangle  \\
&=&(2\pi )^{-n}\left\langle \mathcal{F}(w),M_{\sigma }\circ \mathcal{F}(\varphi )\right\rangle  \\
&=&\left\langle w,\mathcal{F}^{-1}\circ M_{\sigma }\circ\mathcal{F}(\varphi)\right\rangle  \\
&=&\left\langle w,T_{\sigma }\varphi \right\rangle .
\end{eqnarray*}
\qed

\subsection{The Riesz potential and the Riesz operator}

\textbf{Definitions} The \emph{Riesz potential} on $\Phi'$ of
order $\alpha \in \r$ is the operator $I^{\alpha} \in
\mathcal{FM}\Op(\Phi ')$ whose symbol is
$$\Smb(I^{\alpha})= |\xi |^{-\alpha}.$$

\bigskip

\begin{theorem} \label{th.lapinvphi}
$\Delta :\Phi '\rightarrow \Phi '$ is an isomorphism with inverse
$I^{2}:\Phi '\rightarrow \Phi '$.
\end{theorem}

\textbf{Proof} \ We have $\Smb(\Delta )= |\xi |^{2}$, hence 
$\Smb(I^2\circ\Delta) = |\xi |^{-2}\cdot |\xi |^{2} = 1$.

\qed

\bigskip

\begin{corollary} \label{cor.ext.sq1}
$\Delta :S'\rightarrow S'$ is surjective and we we thus have an exact sequence
\begin{equation*}
0\rightarrow \mathcal{H}\rightarrow
\mathcal{S}'\overset{\Delta}{\rightarrow }\mathcal{S}'\rightarrow
0.
\end{equation*}
\end{corollary}

\medskip

\textbf{Proof} The previous theorem says that for any  $f\in
\mathcal{S}'$, we can find a distribution $g\in \mathcal{S}'$ such
that
\begin{equation*}
\Delta g=f\quad \text{in}\quad \Phi '=\mathcal{S}'/
\mathcal{P}.
\end{equation*}
This means that there exists a polynomial $P\in \mathcal{P}$ such that $
\Delta g=f+P$ in $\mathcal{S}'$. By Theorem \ref{Th.InvPoly}, we
can find a polynomial $Q\in \mathcal{P}$ such that $\Delta Q=P$
and it is now clear that
\begin{equation*}
 \Delta (g-Q)=f\quad \text{in}\quad \mathcal{S}'.
\end{equation*}%
This proves that $\Delta (\mathcal{S}')=\mathcal{S}'$.

\qed

\medskip

\textbf{Remark} The distribution $g$ in the above  reasoning  is
only well defined in $\Phi '$ (by the formula $g=I^{2}f$). In the
space $\mathcal{S}'$\ it is only well defined up to a polynomial and
we have no constructive inverse map $\Delta^{-1}:S'\rightarrow S'$.

\medskip

\textbf{Definition} The \emph{Riesz operator} in direction $j$ is the operator
$\mathcal{R}_j \in \mathcal{FM}\Op(\Phi ')$ defined by
$$
 \mathcal{R}_j : = - I^1\circ\partial _{j} = -\partial _{j} \circ I^1.
$$
Its symbol is
$$
 \Smb(\mathcal{R}_j) = -\Smb(I^1)\Smb(\partial _{j})  =
 i  \frac{\xi_j}{|\xi|}.
$$

\medskip

\begin{proposition}\label{prop.riesz}
The Riesz potential and the Riesz operator enjoy the following
properties:
\begin{enumerate}[i)]
\item $I^0 = Id$;
\item $I^{\alpha}\circ I^{\beta } = I^{\beta }\circ I^{\alpha }= I^{\alpha +\beta }$;
\item $I^{-2} =  \Delta =  \sum_j \partial_j^2$;
\item $\Delta \circ I^{\alpha } =I^{\alpha} \circ \Delta = I^{\alpha
-2}$;
\item $\mathcal{R}_i\circ \mathcal{R}_j = \mathcal{R}_j\circ
\mathcal{R}_i = I^2 \partial _{i} \partial _{j}$;
\item $\sum_{j=1}^n\mathcal{R}_j^2 =- Id$;
\item $\langle I^{\alpha }\varphi ,\eta \rangle =
\langle\varphi ,I^{\alpha }\eta \rangle$.
\end{enumerate}
\end{proposition}

The proof is straightforward.

\qed

\medskip

The Riesz potential $I^{\alpha} \in \mathcal{FM}\Op(\Phi ')$ is
sometimes denoted by $I^{\alpha} = \Delta^{-\alpha/2}$, the previous lemma justifies this notation.

\subsection{Convolution operators in $\mathcal{S}$}

Let $T=T_{\sigma} \in  \mathcal{FM}\Op(\Phi ')$ be an operator
such that $\sigma = \Smb (T) \in \mathcal{S}'\cap \mathcal{M}$, then we can
define another operator $\widetilde{T} : \mathcal{S} \to
\mathcal{S}'$ by the convolution
$$\widetilde{T} \varphi = (\mathcal{F}^{-1}{\sigma}) * \varphi.$$
The next lemma is easy to check.
\begin{lemma}
  The relation between $T$ and $\widetilde{T}$ is given by
  $${T}(\varphi) = \widetilde{T}(\varphi) \quad \text{ (mod $\mathcal{P}$)} \
  $$
for any $\varphi \in \Phi$. In other words, the following diagram commutes:
\begin{eqnarray*}
   \mathcal{S}  & \overset{\widetilde{T}}{\longrightarrow} & \mathcal{S}'   \\
  \cup &  & \downarrow \\
  \Phi & \overset{T}{\longrightarrow} & \Phi'
  \end{eqnarray*}
\end{lemma}
\qed

\begin{theorem}
  The symbol of the Riesz potential $I^{\alpha}$ of order $\alpha$
  belongs to $\mathcal{S}'$ if $\alpha < n$. \\
  If $0<\alpha < n$, then $I^{\alpha}$ defines a convolution
  operator $\mathcal{S} \to \mathcal{S}'$ by
  $$I^{\alpha}\varphi = k_{\alpha} * \varphi,$$
  where $k_{\alpha}$ is the Riesz Kernel
  \begin{equation*}
  k_{\alpha}(x) = \frac{1}{\gamma (n,\alpha)} |x|^{\alpha - n}.
\end{equation*}
\end{theorem}

The proof is given in the appendix.

\qed

\bigskip

 The  Riesz potential  
 $I^{\alpha} : \cs \to \cs'$ is  thus given by the explicit  formula
\begin{equation}\label{form.expln}
  I^{\alpha} \varphi (x) =  \frac{1}{\gamma (n,\alpha)}
  \int_{\r^n} \frac{\varphi (y)}{|x-y|^{n-\alpha}} dy
\end{equation}
if $0< \alpha < n$, and 
\begin{equation}\label{form.expl2}
  I^{\alpha} \varphi (x) =  \frac{1}{\gamma (n,\alpha)}
  \int_{\r^n} \varphi (y)\log\frac{1}{|x|} dy
\end{equation}
if $\alpha = n$.

\section{The Lizorkin space and harmonic analysis in $L^{p}(\mathbb{R}^{n})$}

\begin{proposition} The subspace $\Phi \subset L^{p}(\mathbb{R}^{n})$ is  dense
 $1<p<\infty $.
\end{proposition}

The proof  can be found  in \cite[Theorem 2.7]{samko02}.

\bigskip

\begin{proposition}
The space  $L^p(\r^n)$ injects in $\Phi'$ for any $1\leq p < \infty$. More generally, 
 $L^p(\r^n)+L^q(\r^n)$ injects in $\Phi'$ if $p+q < \infty$.
\end{proposition}

\medskip

This result is a direct consequence of the following Lemma:

\begin{lemma}  Let $f,g\in L^{1}_{loc}(\r^n) \cap \mathcal{S}'$ be two locally integrable functions such that 
$$
 \lambda^n \{x\in \r^n \tq  |f(x) - g(x)|  \geq a \} < \infty,
$$
for some $a\geq 0$ (here $\lambda^n$ is the Lebesgue measure). Assume that $f$ and $g$ coincide in $\Phi'$, i.e. 
$$ \int_{\r^n} f \varphi \, dx = \int_{\r^n} g\varphi \, dx $$
for any $\varphi \in \Phi$. Then $(f-g)$ is almost everywhere constant in $\r^n$.
\end{lemma}

\textbf{Proof} This is Lemma 3.8 in \cite{rubin96}. We repeat the proof, which is very short. 
Since $f$ and $g$ coincide in $\Phi'$, we have $P = (f-g) \in  \mathcal{P}$; by hypothesis, we have
$
 \lambda^n \{x\in \r^n \tq  |P(x)|  \geq a \} < \infty,
$
this is only possible if  $P = c$ \ is a constant such that $|c| < a$.

\qed

\medskip

\textbf{Remark.}  The argument also shows that $L^{\infty}(\r^n)/\mathbb{R}$ injects in $\Phi'$.

\bigskip

We summarize the known inclusions in the following lemma:

\begin{lemma}
We have the following inclusions ($1<p<\infty $)%
\begin{equation*}
\Phi \subset \mathcal{S}\subset L^{p}(\mathbb{R}^{n})\subset \Phi '=%
\mathcal{S}'/\mathcal{P}\text{.}
\end{equation*}%
Furthermore $\Phi $ is dense in $L^{p}$ (for the $L^{p}$ norm)
and in $\Phi '$ (for the weak topology).
\end{lemma}

\qed

\begin{lemma}
If $0< \alpha < n$, then the Riesz Kernel $k_{\alpha}$ is a tempered distribution. In fact 
$$ k_{\alpha} \in  \left( L^{r}(\r^n)+ L^{s}(\r^n)\right)  \subset \mathcal{S}'$$
for any $r,s\geq 1$ such that \ $\ds 0 < \frac{1}{s} <  1 - \frac{\alpha}{n} < \frac{1}{r} \leq  1$.
\end{lemma}

\medskip 
\pf Let $\chi_B$ be the characteristic function of the unit ball and set 
 $K_{1} = \chi_B\,  k_{\alpha}$ and
$K_{2} = (k_{\alpha} - K_{1} )=(1-\chi_B)\, k_{\alpha}(x)$.
 It is easy to check that $K_1 \in L^{r}(\r^n)$ for
any $1 \leq r < \frac{n}{n-\alpha}$ and  $K_2 \in L^{s}(\r^n)$ for
any $\frac{n}{n-\alpha} < s < \infty$, thus
$$k_{\alpha}= K_{1}  + K_{2} \in  L^{r}(\r^n)+ L^{s}(\r^n).$$

\qed

\bigskip

\begin{corollary} 
If $0< \alpha < n$, then the Riesz potential  
 \begin{equation*} 
  I^{\alpha} \varphi (x) = k_{\alpha}*\varphi (x)   = \frac{1}{\gamma (n,\alpha)}
  \int_{\r^n} \frac{\varphi (y)}{|x-y|^{n-\alpha}} dy
\end{equation*}
defines a bounded operator 
$$I^{\alpha} : L^{p}(\r^n) \to \left(L^{q_1}(\r^n)+L^{q_2}(\r^n)\right)$$
for any $1\leq p \leq \infty$ and   $1\leq q_1< q_2 < \infty$ such that $\frac{1}{q_2} < \frac{1}{p}  - \frac{\alpha}{n} < \frac{1}{q_1}$.
\end{corollary}

\medskip

In particular, $I^{\alpha}$ is a continuous operator $I^{\alpha} : L^{p}(\r^n) \to \Phi'$ for any $0< \alpha < n$.

\medskip
 
\pf  Let us set $r = (pq_1 - p - q)/(pq_1)$ and  $s = (pq_s - p - q)/(pq_2)$, we then have
\begin{equation}\label{rsyoung}
\frac{1}{r} =  1 + \frac{1}{q_1} - \frac{1}{p}  >  1 - \frac{\alpha}{n}  \quad
\text{ and } \quad  \frac{1}{s} =  1 + \frac{1}{q_2} - \frac{1}{p}  <  1 - \frac{\alpha}{n}
\end{equation}
 
By the Previous Lemma, we may write $k_{\alpha}= K_{1}  + K_{2} $ with
$K_1 \in L^{r}(\r^n)$ and $K_2 \in L^{s}(\r^n)$.

The Young inequality for convolutions says that under the condition (\ref{rsyoung}), we have
\begin{equation*}
 \|K_1*f\|_{L^{q_1}} \leq \|K_1\|_{L^{r}}\|f\|_{L^{p}}
 \quad \text{ and } \quad  
 \|K_2*f\|_{L^{q_2}} \leq \|K_2\|_{L^{s}}\| f\|_{L^{p}}
\end{equation*}
Since $k_{\alpha} = K_{1}  + K_{2} $, we conclude that 
\begin{equation*}
 \|k_{\alpha}*f\|_{L^{q_1}} \leq \left( \|K_1\|_{L^{r}} +
 \|K_2\|_{L^{s}} \right)  \| f\|_{L^{p}}.
\end{equation*}

\qed

\bigskip

If $p \in (1, n/\alpha )$, then we have the following much deeper result:

\begin{theorem}[Hardy-Littlewood-Sobolev]\label{Th.HLS}
The Riesz potential defines a bounded operator $$I^{\alpha} :
L^p(\r^n) \to L^q(\r^n)$$ if and only if $\alpha \in (0,n)$,
$1<p<n/\alpha$ and $q=\frac{np}{n-p\alpha}$. The formulas
(\ref{form.expln}) and (\ref{form.expl2}) still hold in this case.
\end{theorem}

References for this important result are \cite[page 119]{stein70} ,
\cite[Th. 2.2]{samko02}  or \cite[Th. 3.14]{hedberg96}.

\bigskip

Recall that we defined the  Riesz transform in direction $j$ to be the operator 
$ \mathcal{R}_j  = - I^1\circ \partial_j$. 
Its symbol is
$ \rho_j =  i |\xi |^{-1}\xi _{j}$, and  for any $\varphi \in \mathcal{S}$, we thus  have 
$$
 \mathcal{F}(\mathcal{R}_j (\varphi)) = 
   \rho_j  \mathcal{F}(\varphi).
$$
The Riesz transform of a Lizorkin distribution   $f\in \Phi$ is
characterized by 
$$
 \langle \mathcal{R}_j (\varphi) (f), \varphi \rangle =
  \langle f ,   
  \mathcal{F}^{-1}\left(  \rho_j \mathcal{F}(\varphi)\right)
  \rangle 
$$
for any $\varphi \in \mathcal{S}$.

\medskip

The function $\rho_j$ does not belong to the Schwartz space $\mathcal{S}$
and thus its (inverse) Fourier transform  $ \check{\rho}_j = \mathcal{F}^{-1}\left(  \rho_j \right)$ is not a priori  well defined. We can therefore not write  the Riesz transform as a convolution.
However, $ \mathcal{R}_j $ can be represented as a singular integral:

\begin{theorem}[Calderon-Zygmund-Cotlar]\label{Th.CZ}
The Riesz transform $\mathcal{R}_{j}:\mathcal{S}\rightarrow \mathcal{S}'$ is given  by the formula
\begin{equation}\label{RT.int}
\mathcal{R}_{j} \, (\varphi(x)) = \lim_{\delta \rightarrow 0} \ 
\frac{\Gamma(  \frac{n+1}{2} )}{\pi ^{(n+1)/2}}\int_{|y|>\delta}
\frac{(x_{j}-y_{j})}{|x-y|} \varphi(y)dy.
\end{equation}

\medskip

Furthermore, $\mathcal{R}_{j}$ extends as a bounded operator 
$$ \mathcal{R}_j : L^p(\r^n) \to L^p(\r^n)$$
 for all $1<p<\infty$ and the formula
(\ref{RT.int})  still holds in this case.
\end{theorem}

\bigskip

This deep result is a consequence of  sections  II \S 4.2, III \S 1.2 and   III \S 3.3 in  the book of Stein \cite{stein70}, see also  \cite{hedberg96}.

\qed

\medskip

Let us denote by $c_p =  \|\mathcal{R}_{j}\|_{L^p\to L^p}$
the norm of the operator $ \mathcal{R}_j : L^p(\r^n) \to L^p(\r^n)$, it is clearly independant of $j$. 
The exact value of $c_p$ is known, see  \cite[page 304]{iwaniecmartin}; let us only stress that 
$$\lim_{p\to 1} c_p = \lim_{p\to \infty} c_p = \infty.$$

\medskip

\textbf{Remark.} For $p=1$ and $p=\infty$, the Riesz transform is still a
bounded operator in appropriate function spaces, namely
\begin{eqnarray*}
 \mathcal{R}_j : L^1(\r^n) &\to & \text{\small weak}L^1(\r^n)\\
 \mathcal{R}_j : L^{\infty}(\r^n) &\to & BMO(\r^n)
\end{eqnarray*}
are bounded operators.
There are also results on  weighted  $L^p$ spaces satisfying a Muckenhoupt condition.

\subsection{Applications of these $L^p$ bounds} 

To illustrate the power of the two previous theorems, 
we give below  very short proofs  of two important results for
functions in $\r^n$ (compare  \cite[pages 59 and 126]{stein70}).

\begin{theorem}[Sobolev-Gagliardo-Nirenberg]  \label{th.sgn1}
Let $1< p , q < \infty$ be such that  $\frac{1}{p} - \frac{1}{q} = \frac{1}{n}$. There exists a finite constant $C< \infty$ such that  for any $f \in L^p(\r^n)$, 
we have
\begin{equation} \label{ineq.sgn_scal}
  \|f \|_{L^q(\r^n)} \leq C \, \sum_{j=1}^n
 \|\partial_jf \|_{L^p(\r^n)}.
  \end{equation}
\end{theorem}

\medskip

\textbf{Remarks:} 1.) This inequality also holds for $p=1$, see \cite[ chap. V \S 2.5
pp. 128--130]{stein70}, but not for $p=\infty$.

2.) A homogeneity argument shows that the inequality (\ref{ineq.sgn_scal}) cannot hold with a finite constant if  $\frac{1}{p} - \frac{1}{q} \neq \frac{1}{n}$ (see the argument in the proof of Theorem \ref{th.sgn}
below).

\medskip

\pf  Combining the identity 
$ -Id = \sum_j \mathcal{R}_j^2   =  \sum_j I^1 \circ \mathcal{R}_j\circ \partial_j  $
with Theorems \ref{Th.HLS} and \ref{Th.CZ}, we obtain
$$
  \|f \|_{L^q(\r^n)}  \leq   a_p \sum_{j=1}^n
 \| \mathcal{R}_j\partial_jf \|_{L^p(\r^n)}
 \leq   
  a_p c_p\sum_{j=1}^n
 \|  \partial_jf \|_{L^p(\r^n)}.
$$

\qed

\begin{theorem}[A priori estimates for the Laplacian]  \label{th.appest}
The following inequality holds for any $f\in \Phi'$ and any $\mu, \nu = 1,...,n$:
\begin{equation} 
  \|\partial_{\mu} \circ \partial_{\nu}  f \|_{L^p(\r^n)} \leq c_p^2 \, 
 \|\Delta f \|_{L^p(\r^n)}.
 \end{equation}
\end{theorem}

Recall that $c_p < 0$ if and only if $1<p<\infty$.

\medskip

\pf This result is an obvious consequence of  the definition of $c_p$
and the identity
$$\partial_{\mu} \circ \partial_{\nu}  = - \mathcal{R}_{\mu}\circ  \mathcal{R}_{\nu} \circ \Delta.$$
\qed

\bigskip

\textbf{Remark.}  This result holds for any $f \in L^s(\r^n), 1\leq s \leq \infty$, since $L^s(\r^n) \subset \Phi' $. But it does not hold for arbitrary functions 
$f\in \mathcal{S}'(\r^n)$, for instance the harmonic polynomial $f(x,y) = xy
\in \mathcal{H}$\ satisfies $\Delta f = 0$, but $ \partial_x \partial_y f =1 \neq 0$.


\section{Applications to differential forms}

\subsection{Differential forms in $\r^n$}

We denote by $\Lambda^k =  \Lambda^k(\r^{n*})$ the vector space of
antisymmetric multinear $k$-forms on $\r^n$. Recall that
$\dim (\Lambda^k ) = \binom {n}{k}$ and  a basis of  this space is given by 
$$
\{ dx_{i_1}\wedge dx_{i_2} \wedge\cdots\wedge dx_{i_k}  \tq 
{i_1} <  {i_2} < \cdots < {i_k} \}.
$$
A \emph{smooth differential form $\theta$ of degree} $k$ on $\r^n$ is simply a smooth function on $\r^n$ with values in $\Lambda^k$. It is thus uniquely represented as
\begin{equation}\label{form.dif}
 \theta = \sum_{{i_1} <  {i_2} < \cdots < {i_k}}
 a_{{i_1} \cdots  {i_k}}(x)  dx_{i_1}\wedge dx_{i_2} \wedge\cdots\wedge  dx_{i_k}
\end{equation}
 where the coefficients $a_{{i_1} \cdots  {i_k}}$ are smooth functions. We denote by 
$C^{\infty}(\r^n,\Lambda^{k})$ the space of smooth differential forms  of degree $k$ on $\r^n$.

\bigskip

We will also consider later other spaces of differential forms on $\r^n$ such as $L^p(\r^n,\Lambda^{k})$ or $\mathcal{S}(\r^n,\Lambda^{k})$. The form (\ref{form.dif}) belongs to $\mathcal{S}(\r^n,\Lambda^{k})$ if its coefficients $a_{{i_1} \cdots  {i_k}}$
are rapidly decreasing functions and $\theta \in L^p(\r^n,\Lambda^{k})$ if all 
$a_{{i_1} \cdots  {i_k}} \in L^p(\r^n)$.

\bigskip

We shall study a number of operators on differential forms. Observe first that
the operators $\partial_{\mu} = \frac{\partial}{\partial x_{\mu}}$, $I^{\alpha}$ and
$\mathcal{R}_{\mu}$ are well defined on appropriate classes of differential forms
by acting on the coefficients $a_{{i_1} \cdots  {i_k}}$ of the form (\ref{form.dif}).

\bigskip

The \emph{Hodge star operator} is the linear map $\star : \Lambda^k \to \Lambda^{n-k}$ defined by the condition
$$
 \left(dx_{i_1} \wedge\cdots\wedge  dx_{i_k}\right)  \wedge \star\,  \left(dx_{i_1}\wedge   \cdots\wedge dx_{i_k}\right) 
  =   dx_{1}\wedge dx_{2} \wedge\cdots\wedge  dx_{n}
$$
for any  ${i_1} <  {i_2} < \cdots{i_k}$; observe that
\begin{equation}\label{**}
 \star \star = (-1)^{k(n-k)} \func{Id} \qquad \text{on}\  \Lambda^k.
\end{equation}
The $\star$ operator naturally extends to the space
of differential forms with any kind of coefficients.
 
 \bigskip

The \emph{interior product} of the $k$-form $\theta$ with the vector $X$ is the $(k-1)$-form defined by
\begin{equation*}
\interior _{X}\theta (v_{1},...,v_{k-1})=\theta (X,v_{1},...,v_{k-1}).
\end{equation*}%

We  denote by $\interior_{\mu} =\interior_{\frac{\partial}{\partial x_{\mu}}}$  the interior product with $\frac{\partial}{\partial x_{\mu}}$ and by
$\ee_{\mu}$ the exterior product with $dx_{\mu}$:
$$\ee_{\mu} \, \theta = dx_{\mu}\wedge \theta.$$

\begin{lemma}\label{lem.*iota}
The following holds on $k$-forms:
\begin{equation}\label{id.*iota1}
\interior_{\mu} = (-1)^{nk+n} \star \, \ee_{\mu} \star.
\end{equation}%
\end{lemma}

\pf  We first show that for an arbitrary differential form 
$\alpha$, we have
\begin{equation}\label{id.*iota}
 \interior_{\mu} (\star \alpha) =  \star ( \alpha\wedge dx_{\mu}).
\end{equation}%
It is enough to prove this identity  for 
$\alpha =  dx_{j_1} \wedge \cdots \wedge dx_{j_k}$. Observe that
if $\mu = j_r$ for some $r\in \{1,2,\dots, k \}$, then both sides of the 
equation (\ref{id.*iota}) trivially vanish: we thus 
assume $\mu \neq j_r $ for all $r$ and 
set $\beta = \star (\alpha \wedge dx_{\mu})$. Then, by definition
$$\alpha \wedge dx_{\mu} \wedge \beta
=  dx_1 \wedge dx_{2} \wedge \cdots \wedge dx_{n},$$
this relation clearly implies 
$\star \alpha =  dx_{\mu} \wedge \beta$, 
and the equation (\ref{id.*iota}) is now easy to check:
$$
\interior_{\mu} (\star \alpha) = \interior_{\mu} (dx_{\mu} \wedge \beta)
 = \beta = \star (\alpha \wedge dx_{\mu}).
$$
Let us now consider an arbitrary $k-$form $\theta$, and let 
$\alpha =  (-1)^{k(n-k)}\star \theta$, i.e. \ $\star \alpha =
 \theta$. Using  (\ref{id.*iota}), we have 
\begin{eqnarray*}
\interior_{\mu} (\theta) &=& \interior_{\mu} (\star\alpha) =
 \star (\alpha \wedge dx_{\mu})
 =  (-1)^{n-k}  \star (dx_{\mu}\wedge \alpha)
 \\  & = &    (-1)^{n-k}  \star (\ee_{\mu} \, \alpha)
= 
 (-1)^{n-k}(-1)^{k(n-k)}  \star (\ee_{\mu} \star \theta)
  \\  & = &    (-1)^{kn+n}  \star (\ee_{\mu} \star \theta).
\end{eqnarray*}

\qed

\bigskip

We now define the  \emph{exterior differential operator} by 
\begin{equation} \label{id.d}
  d  =  \sum_{\mu=1}^{n}\ee _{\mu} \circ \partial_{\mu}  
   = \sum_{\mu=1}^{n}  \partial_{\mu}  \circ  \ee _{\mu},
   \end{equation}%
 and the  \emph{codifferential operator} by
\begin{equation} \label{id.delta}
  \delta   = - \sum_{\mu=1}^{n}\interior _{\mu} \circ \partial_{\mu}
     = - \sum_{\mu=1}^{n}  \partial_{\mu}  \circ  \interior _{\mu}.
\end{equation}%
It follows from Lemma \ref{lem.*iota}, that for any  $k$-form $\theta$,
we have  
\begin{equation}\label{def.delta}
 \delta\theta =  (-1)^{nk+n+1} \star d\star \, \theta.
\end{equation}
If $\theta$ has the representation (\ref{form.dif})  then $d\theta$  is given by
\begin{equation*}
 d\theta = \sum_{{i_1} <  \cdots < {i_k}}
 da_{{i_1} \cdots {i_k}} \wedge dx_{i_1}\wedge dx_{i_2}\wedge \cdots dx_{i_k},
\end{equation*}
and $\delta\theta$ is given by
\begin{equation*}
 \delta\theta =   \sum_{{i_1} <  \cdots < {i_k}} \sum_{j=1} ^{k} (-1)^j
 \frac{\partial a_{{i_1} \cdots {i_k}}  }{\partial x_{i_j}} \ dx_{i_1}\wedge 
\cdots \wedge \widehat{dx_{i_j}}  \wedge \cdots  \wedge dx_{i_k};
\end{equation*}
A direct  computation show that these operators enjoy the following properties :
\begin{equation*}
d\circ d=\delta \circ \delta =0
\end{equation*}%
and
\begin{equation*}
 \Delta  = (d+\delta )^{2} = d\circ \delta +\delta \circ d = -\sum_{\mu=1}^{n}
 \partial_{\mu}^2
\end{equation*}

\subsection{Temperate and Lizorkin currents}

\textbf{Definitions}  A \emph{rapidly decreasing} differential form of degree $k$ is an element  $\theta \in \mathcal{S}(\r^n, \Lambda^k)$, i.e. a differential form with coefficients in the Schwartz space $\mathcal{S}$.

\medskip

A \emph{temperate current}  $f$ of degree $k$ is a continuous linear form on  $\mathcal{S}(\r^n, \Lambda^k)$, see \cite{schwartz98}. The evaluation of the temperate current $f$ on the 
differential form $\theta \in \mathcal{S}(\r^n, \Lambda^k)$ is denoted by
$$
 \langle f, \theta \rangle \in \mathbb{C}.
$$
the continuity of $f$ means that if $\left\{ \theta_i \right\}\subset \mathcal{S}(\r^n, \Lambda^k)$ is a sequence of rapidly decreasing differential forms such that 
$\theta_i$ converges to $ \theta \in \mathcal{S}(\r^n, \Lambda^k)$ (i.e. all coefficients
converge in the Schwartz space $\mathcal{S}$), then 
$$
 \langle f, \theta \rangle = \lim_{i\to \infty}  \langle f, \theta_i \rangle. 
$$
The space of temperate currents is denoted by $\mathcal{S}' (\r^n, \Lambda^k)$.
Any differential form $f \in L^p(\r^n, \Lambda^k)$ determines a temperate current by
the formula
$$
  \langle f, \theta \rangle = \int_{\r^n} f\wedge * \theta.
$$
This formula defines an embedding $L^p(\r^n, \Lambda^k) \subset \mathcal{S}' (\r^n, \Lambda^k)$.

\bigskip

Another important class of temperate currents is given by the space 
$\mathcal{P}(\r^n, \Lambda^k)$ of differential forms with polynomial coefficients. We may thus define the space of \emph{Lizorkin forms} as 
$$
 \Phi(\r^n, \Lambda^k) = \{\phi \in \mathcal{S}(\r^n, \Lambda^k)
 \tq    \langle P, \phi \rangle = 0 \text{ for all }   P \in \mathcal{P}(\r^n, \Lambda^k)
 \}.
$$
The dual space is called the space of  \emph{Lizorkin currents}, it 
coincides with the quotient
$$
 \Phi'(\r^n, \Lambda^k) = \mathcal{S}' (\r^n, \Lambda^k)/\mathcal{P}(\r^n, \Lambda^k),
$$
we can think of a Lizorkin current as a   differential forms with coefficients in $\Phi'$.

\medskip

Given a Lizorkin current $f \in  \Phi'(\r^n, \Lambda^k)$, we define its
differential $df$, its codifferential $\delta f$, its Laplacian and  its Riesz potential of order $\alpha$ by the following formulas
$$
\begin{tabular}{cc}
$\langle df , \varphi \rangle  =  \langle f , \delta \varphi \rangle$,
& $\langle \delta f , \varphi \rangle  =  \langle f , d\varphi \rangle, $
\\ \\
$\langle \Delta f , \varphi \rangle  =  \langle f , \Delta \varphi \rangle$,  
& \quad  $ \langle I^{\alpha} f , \varphi \rangle  =  \langle f , I^{\alpha} \varphi \rangle$
\end{tabular}
$$
for all $\varphi  \in  \Phi(\r^n, \Lambda^k)$. These are continuous operators
$\Phi'(\r^n, \Lambda^k) \to \Phi'(\r^n, \Lambda^k)$.

\medskip

Observe that  the Riesz potential commutes with $\delta $ and $d$:
\begin{equation*}
\delta I^{2}=I^{2}\delta \qquad \text{and}\qquad dI^{2}=I^{2}d,
\end{equation*}%
and that we have \ $ \Delta I^{2} =I^{2}\Delta = Id$. In particular, we have the
\begin{theorem} 
The Laplacian 
$\Delta : \Phi'(\r^n, \Lambda^k) \to \Phi'(\r^n, \Lambda^k)$ is an isomorphism with inverse
$I^{2}$.
\end{theorem}

\textbf{Proof }  This follows immediately from Theorem \ref{th.lapinvphi}.

\qed

\bigskip
 
\begin{corollary}\label{cor.AA}
We have the exact sequences
\begin{equation*}
0\rightarrow \mathcal{H}(\r^n, \Lambda^k) \rightarrow \mathcal{P}(\r^n, \Lambda^k)\overset{\Delta}{\rightarrow }\mathcal{P}(\r^n, \Lambda^k)\rightarrow 0,
\end{equation*}%
and
\begin{equation*}
0\rightarrow \mathcal{H}(\r^n, \Lambda^k) \rightarrow \mathcal{S}'(\r^n, \Lambda^k)\overset{\Delta}{\rightarrow }\mathcal{S}'(\r^n, \Lambda^k)\rightarrow 0.
\end{equation*}%
In particular, any temperate current  $\theta\in\mathcal{S}'(\r^n, \Lambda^k)$ 
is the sum of an exact  plus a coexact current, more precisely we have
$$\theta = \Delta \omega = \delta(d\omega) + d(\delta\omega) $$
for some $\omega \in\mathcal{S}'(\r^n, \Lambda^k)$,  well defined up to
a harmonic current.
\end{corollary}

\textbf{Proof } The first  exact sequence follows from Theorem \ref{Th.InvPoly} and
the second one is proved as in Corollary  \ref{cor.ext.sq1}, using the previous theorem.

\qed

\bigskip

\textbf{Remark} Observe that the situation is very different from the $L^2$-theory (or any  other Hilbert space model); in $L^2$, we have 
$\overline{\Imm \Delta} = \left(\ker \Delta\right)^{\bot}$; in particular 
$\ker \Delta = 0$ if $\Delta$ is onto, but in $\mathcal{S}'(\r^n, \Lambda^k)$, \
the Laplacian is onto while 
$\ker \Delta = \mathcal{H}(\r^n, \Lambda^k) \neq 0$.

\subsection{The Riesz transform on currents}

\textbf{Definitions.}  We define the \emph{Riesz transform} on 
$\Phi'(\r^n, \Lambda^k)$ by
\begin{equation}\label{expansion.R}
\mathcal{R}  = d\circ I^{1} =  I^{1}\circ d =\sum_{\mu=1}^{n}  \ee_{\mu}\circ
\mathcal{R}_{\mu}
\end{equation}%
and its adjoint
\begin{equation}\label{expansion.R*}
\mathcal{R}^{\ast } =  \delta \circ I^{1} =I^{1}\circ \delta =  - \sum_{\mu=1}^{n}\interior _{\mu}
\circ \mathcal{R}_{\mu}.
\end{equation}
We also define four additional operators 
$E,E^*,U,U^* : \Phi'(\r^n, \Lambda^k) \to \Phi'(\r^n, \Lambda^k)$ by 
\begin{equation}\label{definition.E}
 E = d  \circ \delta \circ  I^2  = \mathcal{R}  \circ \mathcal{R}^*   \qquad \text{,}  \qquad 
 E^* =  \delta \circ d\circ  I^2  =  \mathcal{R}^*  \circ \mathcal{R},
\end{equation}
and
\begin{equation}\label{definition.U}
  U = I^1\circ \mathcal{R}^* = I^2 \circ \delta
 \qquad \text{,}  \qquad 
 U^* = I^1\circ \mathcal{R} = I^2 \circ d.
\end{equation}

\bigskip

\begin{proposition}\label{prop.REU} 
These operators are continuous on $\Phi'(\r^n, \Lambda^k)$. They 
enjoy the following properties:
\begin{enumerate}[(a.)]
\item  $E +E^* =  \mathcal{R}  \circ \mathcal{R}^* + \mathcal{R}^*  \circ \mathcal{R} =  \Delta \circ I^2 = Id$;
\item $E=0$ on $\ker \delta$ and $E^*=0$ on $\ker d$ ;
\item $E\circ E^*=E^* \circ E =0$;
\item $E\circ E=E$ and  $E^*\circ E^*=E^*$;
\item $E=Id$ on $\ker d$ and $E^*=Id$ on $\ker \delta$;
\item $\Imm E  = \ker (E^*) =  \Imm d=  \ker d$ \ and \ $\Imm E^* = \ker E =  \Imm \delta =  \ker \delta$;
\item $(\mathcal{R},\mathcal{R}^*)$ \ and \ $(U,U^*)$ are adjoint pairs , i.e. 
 $$
\langle \mathcal{R}  \theta , \varphi \rangle = \langle  \theta  , \mathcal{R}^{\ast }\varphi  \rangle \quad  \mathrm{and } \quad
\langle \mathcal{R}^{\ast }  \theta  , \varphi \rangle = \langle f \theta , \mathcal{R} \varphi  \rangle,
$$
for any $\theta \in \Phi'(\r^n, \Lambda^k)$ and $\varphi \in \Phi (\r^n, \Lambda^k)$, and likewise for $U$.\\

\item $E$ and $E^*$ are self-adjoint.
\item $E = d \circ U$ and $E^* =  \delta \circ U^*$.
\end{enumerate}
\end{proposition}

\pf   (a) Follows from the definitions and $Id = \Delta \circ I^2 = (d\delta + \delta d)\circ I^2 
= E +E^*$.\\
(b) If $\delta \, \theta = 0$, then $E\, \theta = d\, \delta  \,  I^2 \, \theta
=  I^2  \, d\, \delta \, \theta = 0$, hence  $E^*=0$ on $\ker d$.  A similar argument shows that  
 $E^*=0$ on $\ker d$. \\
 (c) By definition $E\circ E^*= d \, \delta  \, I^2  \, \delta \, d \,   I^2
= I^4 \, d \, \delta ^2 d   = 0$. The proof that $E^* \circ E =0$
is the same. \\
(d) This follows  from (a) and (c), since
$$E = E\circ Id = E\circ E + E \circ  E^* = E\circ E.$$
(e) This follows immediately from (a) and (b). \\
(f) From $E\circ E = E$, we have $\theta \in \Imm E$ if and only if  $\theta = E\theta$. Since 
 $E+E^* = Id$, we have $E^*\theta = (Id-E)\theta = 0$, thus $\Imm E = \ker E^*$. 
 Furthermore, using $E = d\delta I^2$  and  Property (e),   we see that
$$
 \Imm E  \subset  \Imm  d \subset \ker d  = E(\ker d) \subset \Imm E.
$$
This shows that $\Imm E = \Imm d=  \ker d$.
The proof that $\Imm E^* = \ker E = \Imm \delta =  \ker \delta$ is similar. \\
(g) We have  
\begin{equation*}
\langle \mathcal{R} \theta , \varphi \rangle =
\langle d I^1\theta ,\varphi \rangle = 
\langle I^1\theta , \delta\varphi \rangle =
\langle \theta , I^1 \delta \varphi \rangle =
\langle \theta , \mathcal{R}^*\varphi \rangle.
\end{equation*}%
The proof that $(U,U^*)$ is an adjoints pair is similar, using $U=I^2 \delta$. \\
(h) It follows that  $E$ is selfadjoint, for 
\begin{equation*}
\langle E \theta , \varphi \rangle =
\langle \mathcal{R} \mathcal{R}^*\theta ,\varphi \rangle = 
\langle \mathcal{R}^*\theta , \mathcal{R}^*\delta\varphi \rangle =
\langle \theta , \mathcal{R} \mathcal{R}^*\varphi \rangle =
\langle \theta , E\varphi \rangle,
\end{equation*}
and likewise for $E^*$. \\
(i) We have $dU = d I^2 \delta = d \delta I^2 = E$ and $\delta U^* = \delta I^2 d =  \delta d I^2 = E^*$.

\qed

\bigskip
Let us denote by 
$$
  E\Phi'(\r^n, \Lambda^k) = d \Phi'(\r^n, \Lambda^{k-1})
$$
the space of  \emph{exact Lizorkin currents} of degree $k$ and by 
$$
  E^*\Phi'(\r^n, \Lambda^k) = \delta \Phi'(\r^n, \Lambda^{k+1})
$$
the space of  \emph{coexact Lizorkin currents} of degree $k$.

\medskip

\begin{corollary}  These subspaces can be expressed as
   $$
    E\Phi'(\r^n, \Lambda^k) = \Imm (E) =  \ker E^*
     = \ker \left[ d : \Phi'(\r^n, \Lambda^k) \to \Phi'(\r^n, \Lambda^{k+1})\right]
    $$ 
and
   $$
    E^*\Phi'(\r^n, \Lambda^k) = \Imm (E^*) = \ker E  = \ker \left[ \delta : \Phi'(\r^n, \Lambda^k) \to \Phi'(\r^n, \Lambda^{k-1})\right].
    $$
In particular    $E\Phi'(\r^n, \Lambda^k)$ and $E^*\Phi'(\r^n, \Lambda^k)$ are closed subspaces in
   $\Phi'(\r^n, \Lambda^k)$ and  we have a direct sum decomposition
   $$\Phi'(\r^n, \Lambda^k) =   E\Phi'(\r^n, \Lambda^k) \oplus   E^*\Phi'(\r^n, \Lambda^k).$$
 \end{corollary}
\textbf{Proof} This is obvious from the  previous proposition.

 \qed

\medskip

\textbf{Remarks.}  \textbf{(1.) }Thus $E$ and $E^*$  are the projections 
of $\Phi'(\r^n, \Lambda^k)$ onto $E\Phi'(\r^n, \Lambda^k)$ and
$E^*\Phi'(\r^n, \Lambda^k)$ respectively.  One says that $E\theta$ is the \emph{exact part} of $\theta \in \Phi'(\r^n, \Lambda^k)$ and $E^*\theta$ is its \emph{coexact part}. The formula 
$$\theta = E\theta + E^*\theta$$
is the \emph{Hodge-Kodaira decomposition} of the Lizorkin distribution $\theta$.

\medskip

\textbf{(2.)} The last part says in particular that there is no cohomology in $\Phi'(\r^n, \Lambda^k)$, i.e. 
$$
 \cdots \to \Phi'(\r^n, \Lambda^{k-1})  \stackrel{d}{\to} \Phi'(\r^n, \Lambda^{k}) \stackrel{d}{\to} \Phi'(\r^n, \Lambda^{k+1}) \to \cdots
$$
is an exact sequence.

\medskip

\textbf{(3.)}  Using the equalities $U = I^2 \circ \delta$ and $U^* =   I^2 \circ d$, and observing
that  $E = d \circ U$ and  $E^* =  \delta \circ U^*$, one can write  the Hodge-Kodaira decomposition of 
$\theta \in \Phi'(\r^n, \Lambda^k)$ as
$$\theta = d(U\theta) + \delta (U^*\theta).$$

\subsection{Proof of Theorems  \ref{th.hodgeLp} and \ref{th.hodgeLqp}}
\label{sec.LpTheory}

The operators defined in the previous section are well behaved on $L^p$:

\begin{theorem}\label{RLpLk}
The Riesz transform  and its dual
$$\mathcal{R},\mathcal{R}^*:L^{p}(\r^n,\Lambda ^{k})\rightarrow L^{p}(\r^n,\Lambda ^{k-1})$$
are bounded operators on $L^p$ for any $1<p<\infty$.
\end{theorem}

\pf The boundedness of these operators on  $L^p$  follows from Theorem \ref{Th.CZ} and the expansions (\ref{expansion.R}) and (\ref{expansion.R*}). 

\qed  

\begin{theorem}\label{ULpqk}
The operators $U^*$ and $U$  restrict as bounded operators
\begin{equation*}
U^*,U:L^{p}(\r^n,\Lambda ^{k})\rightarrow L^{q}(\r^n,\Lambda ^{k-1}),
\end{equation*}%
if and only if $\frac{1}{p} - \frac{1}{q} = \frac{1}{n}$.
\end{theorem}

\pf This follows from Theorems \ref{Th.HLS}.

\qed  

\bigskip
 
We can now prove the  Theorems stated in the introduction.

\subsection*{Proof of Theorem  \ref{th.hodgeLp}}   
The equation (\ref{eq.HKDLP}) is a trivial consequence of the equality $E+E^* =  \mathrm{Id}$ and
$\Imm (E) \cap \Imm (E^*) = \{ 0 \}$. 
We know from Theorem \ref{RLpLk}, that the  Riesz transform  and its dual are well defined
 bounded operators $\mathcal{R},\mathcal{R}^*:L^{p}(\r^n,\Lambda ^{k})\rightarrow L^{p}(\r^n,\Lambda ^{k-1})$. The operators 
$E = \mathcal{R}\circ \mathcal{R}^*$ and $E^* = \mathcal{R}^*\circ \mathcal{R}$ are then also clearly bounded on  $L^{p}(\r^n,\Lambda ^{k-1})$. 
The algebraic properties  (iii), (iv) and (v) in  Theorem  \ref{th.hodgeLp} are proved in Proposition \ref{prop.REU}, and we know that   $EL^p(\r^n,\Lambda^k)\subset L^p(\r^n,\Lambda^k)$ is a closed subspace,
since it coincides with the kernel of the bounded operator  $E^* : L^p(\r^n,\Lambda^k) \to L^p(\r^n,\Lambda^k)$. Likewise $E^*L^p(\r^n,\Lambda^k)=\ker E$ is also bounded in $L^p(\r^n,\Lambda^k)$.

\qed

\subsection*{Proof of Theorem  \ref{th.hodgeLqp}}
By  Theorem \ref{ULpqk}, the
operators $U^*,U$  restrict as bounded operators on $L^{q}(\r^n,\Lambda ^{k-1})$. The 
relations $dU = E$ and  $\delta U^* = E^*$ are given in Proposition \ref{prop.REU}. 

\qed

\section{Some additional applications}

\subsection{The Gaffney inequality}

\begin{theorem} 
Assume  $1<p <\infty$. There exists a constant $C_p$ such  for any 
 $\theta \in L^{p}(\r^n,\Lambda ^{k})$  and 
any $\mu = 1,2,\cdots n$, we have
\begin{equation} 
  \|\partial  _{\mu}\theta \|_{L^p(\r^n)} \leq  C_p\, \left( 
   \|d\theta \|_{L^p(\r^n)} +  \|\delta \theta \|_{L^p(\r^n)} \right).
\end{equation}
\end{theorem} 

\pf   By  Theorem \ref{Th.CZ} , we know that $\mathcal{R}_{\mu}$ is a bounded operator on $L^p$  and by Theorem  \ref{RLpLk}, it is also   the case of $\mathcal{R}$ and $\mathcal{R}^*$.
The Theorem follows now immediately from  the  following identity:
\begin{equation} \label{id.pregaffney}
  \partial  _{\mu} = \mathcal{R} _{\mu}\circ \mathcal{R} \circ \delta + \mathcal{R} _{\mu}\circ \mathcal{R}^* \circ d.  
\end{equation}

The latter formula is a consequence of the relations
 $\Delta\circ I^2 = Id$ and $\mathcal{R} _{\mu} = I^1\circ \partial  _{\mu}$,
 indeed:
 \begin{eqnarray*}
\partial _{\mu}  &=&  \Delta\circ I^2 \circ \partial _{\mu} 
 \\ &=& 
 \left( I^1 \circ \partial _{\mu} \right)\circ \left(\Delta\circ I^1 \right)
 \\ &=& 
  \mathcal{R} _{\mu} \circ \left(d \delta + \delta d \right)\circ I^1
 \\ &=& 
   \mathcal{R} _{\mu} \circ \left(d\circ I^1 \right)\circ  \delta + 
   \mathcal{R} _{\mu} \circ \left(\delta\circ I^1 \right)\circ d
 \\ &=& 
  \mathcal{R} _{\mu}\circ \mathcal{R} \circ \delta + \mathcal{R} _{\mu}\circ \mathcal{R}^* \circ d.  
\end{eqnarray*}
\qed

\subsection{A Sobolev inequality for differential forms} 

We have the following Sobolev-Gagliardo-Nirenberg inequality for differential forms
on $\r^n$:
\begin{theorem}  \label{th.sgn}
Let $1< p,q< \infty$. There exists a constant $C< \infty$ such that  for any $\theta \in L^p(\r^n,\Lambda ^{k})$, 
we have
\begin{equation} \label{ineq.sgn}
 \|\theta \|_{L^q(\r^n)} \leq C \, \left( 
   \|d\theta \|_{L^p(\r^n)} +  \|\delta \theta \|_{L^p(\r^n)} \right).
\end{equation}
if and only if $\frac{1}{p} - \frac{1}{q} = \frac{1}{n}$.
\end{theorem}

\pf We have $\theta = E\theta + E^*\theta = U(d\theta) + U^*(\delta \theta)$.
By  Theorem \ref{th.hodgeLp}, we know that $U,U^* : L^{p}(\r^n,\Lambda^{k-1}) \to L^{q}(\r^n,\Lambda^{k-1})$ are bounded operators if and only if $\frac{1}{p} - \frac{1}{q} = \frac{1}{n}$, hence the inequality (\ref{ineq.sgn}) holds with
$$
 C =  \max \left\{\|\mathcal{U}\|_{L^p\to L^q}  , \|\mathcal{U}^*\|_{L^p\to L^q} \right\}.
$$
We need to show in the converse direction,that the inequality (\ref{ineq.sgn}) cannot hold with a finite constant if  $\frac{1}{p} - \frac{1}{q} \neq \frac{1}{n}$. 
To do that, we consider  a non zero $k$-form $\theta \in L^p(\r^n,\Lambda ^{k})$,
observe that either $d\theta \neq 0$ or $\delta\theta \neq 0$, for, otherwise
the form $\theta$ would have constant coefficients,  which is impossible for
a  non zero form in  $L^p(\r^n)$. The following quantity is therefore well defined:
$$
 Q(t) = 
    \frac {\|h_t^*\, \theta  \|_{L^q(\r^n)}}{\|h_t^*d\theta \|_{L^p(\r^n)} +  \|h_t^*\delta \theta \|_{L^p(\r^n)}},
$$
where $h_t$ is the $1$-parameter group of linear dilations in $\r^n$ given by $h_t(x) = t \cdot  x$.

A calculation shows that for any $\omega \in L^s(\r^n,\Lambda ^{m})$, we have 
\begin{equation}\label{id..homg1}
   \|h_t^*\, \omega \|_{L^s(\r^n)} = t^{m-\frac{n}{s}} \|\, \theta  \|_{L^s(\r^n)},
\end{equation}
since $d h_t^*\, \theta  = h_t^*d\, \theta$  is a $(k+1)$-form, we obtain
\begin{equation}\label{id..homg2}
  \|h_t^*d\theta \|_{L^p(\r^n)}  
  = t^{1+k-\frac{n}{p}} 
   \|d\theta \|_{L^p(\r^n)}.
\end{equation}
Now be careful, because \ $\delta h_t^* \neq h_t^*\delta$. In fact
 $\delta h_t^*\, \theta  = t^2\, h_t^*\delta \, \theta$; this is a $(k-1)$-form and 
 thus (\ref{id..homg1}) implies that
\begin{equation}\label{id..homg3}
 \|h_t^*\delta \theta \|_{L^p(\r^n)} 
= t^2 t^{k-1-\frac{n}{p}}   \|\delta \theta \|_{L^p(\r^n)} 
= t^{1+k-\frac{n}{p}}   \|\delta \theta \|_{L^p(\r^n)} 
\end{equation}
The last three  identities give us
$$
 Q(t) = \frac{t^{k- \frac{n}{q}}}{t^{1+k-\frac{n}{p}}} \;  Q(1)
 = t^{\frac{n}{p} - \frac{n}{q} - 1} \,  Q(1).
$$
If $\frac{1}{p} - \frac{1}{q} - \frac{1}{n} < 0$, then $\lim_{t \to 0} Q(t) = \infty$
and if $\frac{1}{p} - \frac{1}{q} - \frac{1}{n} > 0$, then $\lim_{t \to \infty} Q(t) = \infty$. We conclude that the Sobolev inequality (\ref{ineq.sgn}) cannot hold if $\frac{1}{p} - \frac{1}{q} - \frac{1}{n} \neq 0$.

\qed

\subsection{The $L^p$ \textit{a priory } estimates for the Laplacian on forms }

\begin{theorem}
The following inequality holds for any  $\theta \in \Phi'(\r^n,\Lambda ^{k})$,
 $1<p <\infty$ and  $\mu,\nu = 1,2,\cdots n$:
\begin{equation} 
  \|\partial  _{\mu} \partial  _{\nu} \theta \|_{L^p(\r^n)} \leq c_p^2 \, 
 \|\Delta \theta \|_{L^p(\r^n)} 
\end{equation}
where  $c_p$ is 
the norm of the operator $\mathcal{R}_j : L^p(\r^n) \to L^p(\r^n)$.
\end{theorem}

Observe that this estimate is actually a scalar estimate.

\medskip

\pf By  Theorem \ref{Th.CZ}, we know that $\mathcal{R}_{\mu}$ is a bounded operator on $L^p$. It is also clearly the case of $\mathcal{R}$ and $\mathcal{R}^*$. The Corollary now follows from the following Lemma and the definition of $c_p$.

\qed

\begin{lemma} 
The following identity holds in $ \Op\Phi'(\r^n,\Lambda ^{k})$:
\begin{equation} \label{ }
  \partial  _{\mu} \partial  _{\nu}  =    \mathcal{R} _{\mu}\circ \mathcal{R} _{\nu}\circ   \Delta.
\end{equation}
\end{lemma}

\pf We have 
\begin{eqnarray*}
\partial  _{\mu} \partial  _{\nu}  &=&  \Delta\circ I^2 \circ \partial  _{\mu} \circ \partial  _{\nu} 
 \\ &=& 
\Delta\circ (I^1 \circ \partial  _{\mu}) \circ (I^1 \circ \partial  _{\nu} )
 \\ &=& 
  \Delta \circ \mathcal{R} _{\mu}\circ \mathcal{R} _{\nu}
  \\ &=& 
  \mathcal{R} _{\mu}\circ \mathcal{R} _{\nu}\circ   \Delta.
\end{eqnarray*}

\qed

\subsection{The $L_{q,p}$-cohomology of $\mathbb{R}^n$}

The set of closed forms in $L^p(\r^n,\Lambda^k)$ is  denoted by 
$$Z^k_p(\r^n) =L^p(\r^n,\Lambda^k) \cap \ker d,$$
and the set of exact forms in $L^p(\r^n,\Lambda^k)$ which are differentials of
forms in $L^q$ is denoted by
\[
  B^k_{q,p}(\r^n) = d\left(L^{q}(\r^n,\Lambda^{k-1}) \right) \cap
  L^p(\r^n,\Lambda^k).
\]
 
\begin{lemma}
  $Z^k_p(\r^n) \subset L^p(\r^n,\Lambda^k)$  is a closed linear subspace.
\end{lemma}
\textbf{Proof.}  By definition, a form $\theta \in L^p(\r^n,\Lambda^k)$ belongs to 
$Z^k_p(\r^n)$ if and only if $\int_{\mathbb{R}^n}\theta\wedge \varphi = 0$ for any
$\varphi \in \mathcal{S}(\r^n,\Lambda^{n-k})$. Suppose now that $\theta \in L^p(\r^n,\Lambda^k)$
 is in the closure of $Z^k_p(\r^n)$. This means that there
exists a sequence $\theta_i \in Z^k_p(\r^n)$ converging to $\theta$ for the $L^p$
norm. Using the H\"older inequality with $q=p/(p-1)$, we have
$$
\left| \int_{\mathbb{R}^n}\theta\wedge \varphi \right| =
 \lim_{i\to \infty} \left| \int_{\mathbb{R}^n}(\theta - \theta_i)\wedge \varphi  \right| \leq  \lim_{i\to \infty} \| \theta - \theta_i\|_{L^p}  \| \varphi\|_{L^q}  = 0,
$$
and therefore $\theta \in Z^k_p(\r^n)$.

\qed

\medskip
 
\textbf{Remark.} It is clear that
$
  B^k_{q,p}(\r^n)  \subset EL^p(\r^n,\Lambda^k) \subset Z^k_p(\r^n).
$
By Proposition \ref{prop.REU} (a), \ref{prop.REU} (b) and Theorem \ref{th.hodgeLp},
we have in fact $Z^k_p(\r^n) = EL^p(\r^n,\Lambda^k)$.  
 
 \bigskip
 
\textbf{Definition.}
The $L_{q,p}$\emph{-cohomology} of $\r^n$ is the quotient
$ \ds
H_{q,p}^{k}(\r^n)= Z_{p}^{k}(\r^n)/B_{q,p}^{k}(\r^n).
$

\medskip

The next result computes this cohomology:  

\begin{theorem}
For any $p,q \in (1,\infty)$ and $1\leq k \leq n$, we have
\[
  H_{q,p}^{k}(\r^n)  = 0 \ \Leftrightarrow  \  \frac{1}{p} - \frac{1}{q} = \frac{1}{n}.
\]
\end{theorem}

\pf 
Assume first that $\frac{1}{p} - \frac{1}{q} = \frac{1}{n}$. By  Proposition \ref{prop.REU}, we have for any $\theta \in Z_{p}^{k}(\r^n)$
$$\theta = E\theta + E^*\theta = E\theta = d(U\theta),$$
because $E^* = 0$ on $\ker d$. By  Theorem \ref{th.hodgeLp}, we know that $U :L^{p}(\r^n,\Lambda^{k-1}) \to L^{q}(\r^n,\Lambda^{k-1})$ is a bounded operator.
Hence $\theta =  d(U\theta) \in B_{q,p}^{k}(\r^n)$; but since  $\theta \in Z_{p}^{k}(\r^n)$ is arbitrary, we have
$$
 H_{q,p}^{k}(\r^n) = Z_{p}^{k}(\r^n)/B_{q,p}^{k}(\r^n) = 0.
$$

To prove the converse direction, we use the interpretation of the $L_{q,p}$-cohomology in terms of Sobolev inequalities.
In particular, it is proven in   \cite[Theorem 6.2]{GT2006} that  if $H_{q,p}^{k}(\r^n) = 0$, then 
 there exists a constant $C$ such that for any $\phi \in  L^q(\r^n,\Lambda^{k-1})$, there  exists a  closed form $\zeta = \zeta(\phi) \in Z_{q}^{k-1}(\r^n)$ such that
 \begin{equation}\label{inconc.sob2}
\left\Vert \phi-\zeta\right\Vert _{L^{q}} \leq C\left\Vert
d\phi\right\Vert _{L^{p}}.
\end{equation}

Let us fix a form $\phi \in  L^q(\r^n,\Lambda^{k-1})$ which is \emph{not} closed and apply the above
inequality to $h_t^*\, \phi$, where $h_t(x) = t \cdot  x$. It says in this case that for any $t\in \r$, there exists 
$\zeta_t\in Z_{q}^{k-1}(\r^n)$ such that
\begin{equation}\label{inconc.sob2a}
\left\Vert  h_t^*\, \phi-\zeta_t\right\Vert _{L^{q}} \leq C\left\Vert
h_t^*\,d\phi\right\Vert _{L^{p}}.
\end{equation}
Using the identity (\ref{id..homg1}) with $s=q,m=(k-1)$ and $s=p, m= k$, we obtain
the inequality
\begin{equation}\label{inconc.sob2b}
\left\Vert  \phi- h_{-t}^*\zeta_t\right\Vert _{L^{q}} \leq C t^{\gamma}\, \left\Vert
d\phi\right\Vert _{L^{p}}
\end{equation}
with $\gamma = 1 + \frac{n}{q} - \frac{n}{p}$. The right hand side of this inequality converges to
zero as $t\to 0$ if $\gamma <0$ or as $t\to \infty$ if $\gamma >0$. Since 
$h_{-t}^*\zeta_t \in  Z_{q}^{k-1}(\r^n)$ for any $t$ and $Z_{q}^{k-1}(\r^n)\subset  L^q(\r^n,\Lambda^{k-1})$
is closed, it follows that $\phi \in Z_{q}^{k-1}(\r^n)$. But $\phi$ is not closed by hypothesis, 
we thus conclude that $\gamma = 0$.
To sum up, the argument shows that if $H_{q,p}^{k}(\r^n) = 0$, then $\gamma = 1 + \frac{n}{q} - \frac{n}{p} = 0$.

\qed

\bigskip

\textbf{Remark.} Theorem 6.2 in \cite{GT2006} says in fact that there exists a constant $C$ such that (\ref{inconc.sob2}) holds
if and only if $T_{q,p}^{k}(\r^n) = 0$ (provided $1< q,p < \infty$). Here $T_{q,p}^{k}$ is the \emph{torsion} which is defined
to be the quotient $\overline{B_{q,p}^{k}}/B_{q,p}^{k}$. We have thus also proved that 
$T_{q,p}^{k}(\r^n)  = 0 \ \Leftrightarrow  \  \frac{1}{p} - \frac{1}{q} = \frac{1}{n}$.
In particular $H_{q,p}^{k}(\r^n)$ is infinite dimensional if $\frac{1}{p} - \frac{1}{q} \neq \frac{1}{n}$.

\section*{Appendix: Computation of the Fourier transform of the Riesz kernel}

\textbf{Definition} The \emph{Riesz kernel} of order $\alpha  \in (0,n)$ is the function
$k_{\alpha}$ defined on $\r^n$ by
\begin{equation*}
  k_{\alpha}(x) = \frac{1}{\gamma (n,\alpha)} |x|^{\alpha - n},
\end{equation*}
where the normalizing constant is given by
$$
 \gamma (n,\alpha) =
 2^{\alpha}\pi^{n/2}\frac{\Gamma (\frac{\alpha}{2})}{\Gamma
 (\frac{n-\alpha}{2})}.
$$

\begin{theorem} \label{four.ka}
The Fourier transform of  the Riesz kernel of order $\alpha  \in (0,n)$ is given by
\begin{equation*}
  \mathcal{F}(k_{\alpha})= |\xi |^{-\alpha }.
\end{equation*}%
\end{theorem}

\medskip

\textbf{Proof}  We will use the fact that the Gaussian function $g(x)=e^{-s|x|^{2}}$
belongs to $\mathcal{S}$ for any $s>0$ and that its Fourier
transform is given by
\begin{equation}\label{ff}
\mathcal{F}(e^{-s|x|^{2}})(\xi )= \left(\frac{\pi
}{s}\right)^{n/2}e^{-|\xi |^{2}/4s} ,
\end{equation}%
(this is a well known fact. see e.g. \cite[Proposition 8.24] {folland99} or
\cite[page 38] {strichartz94}).

\bigskip

To compute the Fourier transform of $k_{\alpha }$, we start from
the formulas
\begin{equation}\label{For.Gammapower}
\Gamma (z)\; a^{-z} = \int_{0}^{\infty}s^{z-1}e^{-as}ds \quad
\mathrm{ and } \quad \Gamma (w)\; b^{-w} = \int_{0}^{\infty
}s^{-w-1}e^{-b/s}ds,
\end{equation}%
which hold for any  $a,b\in (0,\infty)$ and any $z,w\in
\mathbb{C}$ such that $\func{Re}(z), \func{Re}(w)>0$.

\medskip

To check these formulas, use  the substitution $t=as$ (for the
first identity) and   $t=b/s$ (for the second identity) in the
definition $\Gamma (z)=\int_{0}^{\infty }t^{z-1}e^{-t}dt$ of the
Gamma function.

\medskip

We will use the first formula with $a=|x|^{2}$ and
apply the Fourier transform;  keeping in mind the identity
(\ref{ff}), we have
\begin{eqnarray*}
\mathcal{F}\left(\Gamma (z)|x|^{-2z }\right) &=&\mathcal{F}\left(
\int_{0}^{\infty }s^{z-1}e^{-|x|^{2}s}ds\right)  \\
&\overset{Fubini}{=}&\int_{0}^{\infty }s^{z-1}
\mathcal{F}\left( e^{-|x|^{2}s}\right) ds \\
&=&\int_{0}^{\infty }s^{z-1}
\left(\frac{\pi }{s}\right)^{n/2}e^{-|\xi |^{2}/4s} ds \\
&=& \pi ^{n/2} \int_{0}^{\infty }s^{z-\frac{n}{2} -1}
e^{-|\xi|^{2}/4s}ds.
\end{eqnarray*}%
Setting $b=|\xi |^{2}/4$ and $w = \frac{n}{2} - z$, we obtain
from  the second identity in (\ref{For.Gammapower})
\begin{equation*}
  \Gamma (z) \mathcal{F}\left(|x|^{-2z }\right) =
  \pi ^{n/2} \int_{0}^{\infty }s^{-w-1}e^{-b/s}ds
  =  \pi ^{n/2}\Gamma (w)4^w \; |\xi|^{-2w}.
\end{equation*}
Let us set $\alpha = n-2z$, thus $z=\frac{n-\alpha}{2}$ and $w =
\frac{n}{2} - z=\frac{\alpha}{2}$; we write this formula as
\begin{equation*}
\mathcal{F}(|x|^{\alpha -n})=\gamma (n,\alpha )|\xi |^{-\alpha},
\end{equation*}%
where%
\begin{equation*}
\gamma (n,\alpha )=\pi ^{n/2}2^{\alpha }\frac{\Gamma
(\frac{\alpha }{2})}{\Gamma (\frac{n-\alpha }{2})}.
\end{equation*}
The above calculation assumes $\func{Re}(z), \func{Re}(w)>0$, which is equivalent to  $0< \alpha < n$.

 \qed

\bigskip

\textbf{Remark.}  Using the Fourier transform in the Lizorkin sense, it is possible to extend
the Riesz kernel $ k_{\alpha}$ of order $\alpha$ for any real number $\alpha >0$ (and in fact any complex number with $\text{Re}\,  \alpha > 0$). We define it as follow \begin{equation*}
  k_{\alpha} = \frac{1}{\gamma (n,\alpha)} |x|^{\alpha - n}
\end{equation*}
if $\alpha \neq n + 2m$ for any $m\in \mathbb{N}$, and by
\begin{equation*}
  k_{\alpha} = \frac{1}{\gamma (n,\alpha)} |x|^{\alpha - n}\, \log
  \frac{1}{|x|}
\end{equation*}
if $\alpha = n + 2m$ for some $m\in \mathbb{N}$.

\bigskip

With this definition, the previous result is still valid
\begin{proposition}
 The Fourier transform of $k_{\alpha}$, $\alpha \in \mathbb{C}$ is given by
\begin{equation*}
\mathcal{F}(k_{a}) =|\xi |^{-\alpha }.
\end{equation*}
\end{proposition}

\qed


\vfill

Marc Troyanov     \\
Section de Mathmatiques \\
 \'Ecole Polytechnique F{\'e}derale de
Lausanne, 
\\ 1015 Lausanne - Switzerland
\\ marc.troyanov@epfl.ch

\medskip

AMS subjclass:  {58A10, 42B, 42B20, 58A14.}

Keywords: Hodge decomposition, temperate currents, Lizorkin currents.



\begin{thebibliography}{mm}
\bibitem{AR}  S. Axler and W. Ramey,  \textit{Harmonic polynomials and
Dirichlet-type problems}. Proc. Amer. Math. Soc. 123 (1995), no.
12, 3765--3773.

\bibitem{conway} J. Conway  \textit{A course in functional analysis.}  Graduate Texts in Mathematics, \textbf{96}. Springer-Verlag.

\bibitem{derham} G. de Rham \emph{Differentiable manifolds. Forms, currents, harmonic
forms}. Grundlehren der Mathematischen Wissenschaften \textbf{266}.
Springer-Verlag, Berlin, 1984.

\bibitem{folland99} G. B. Folland,
 \textit{Real analysis. Modern techniques and their applications}. Second edition. Pure and Applied Mathematics (New York). A Wiley-Interscience Publication, (1999).
 
\bibitem{GT2006}  V. Gol'dstein and M. Troyanov,
\textit{Sobolev inequalities for differential forms and  $L_{q,p}$-cohomology.}
J. Geom. Anal., Vol. 16, Number 4, 2006.

\bibitem{iwaniecmartin} T. Iwaniec and G, Martin,
\textit{Geometric function theory and non-linear analysis.
Oxford Mathematical Monographs.}
Oxford University Press, New York, 2001.

\bibitem{hedberg96} A. David and R. Hedberg,  
\textit{Function spaces and potential theory}. Grundlehren der Mathematischen Wissenschaften 314, Springer-Verlag, Berlin, 1996.

\bibitem{rubin96} B. Rubin,
\textit{Fractional integrals and potentials}.
Pitman Monographs and Surveys in Pure and Applied Mathematics, 82.
Longman, Harlow, 1996.

\bibitem{samko77}  S. Samko,
\textit{Test functions that vanish on a given set, and division by a function}. (Russian) Mat. Zametki 21 (1977), no. 5, 677--689
(English translation: Math. Notes 21 (1977), no. 5--6, 379--386.)

\bibitem{samko02} S. Samko, 
\textit{Hypersingular integrals and their applications.}
Analytical Methods and Special Functions, 5.
Taylor \& Francis, Ltd., London, 2002.

\bibitem{scarfiello} R. Scarfiello,
\textit{Sur le changement de variables dans les distributions et leurs transformes de Fourier.} 
Nuovo Cimento, IX. Ser. 12, 471-482 (1954). 
\bibitem{scott95} C. Scott, 
\emph{$L\sp p$ theory of differential forms on manifolds.}  
Trans. Amer. Math. Soc. 347 (1995), no. 6, 2075--2096. 
\bibitem{schwartz57} L. Schwartz, 
\textit{Thorie des distributions  valeurs vectorielles. I.}
Ann. Inst. Fourier, Grenoble 7 1957 1--141.

\bibitem{schwartz98} L. Schwartz, 
\textit{Thorie des distributions.}  Nouveau tirage. 
 Paris: Hermann. (1998).
 
 \bibitem{schwarz95} G. Schwarz,
\textit{Hodge decomposition---a method for solving boundary value problems.}
Lecture Notes in Mathematics, 1607,
Springer-Verlag, Berlin, 1995.

\bibitem{stein70} E. Stein,  \textit{Singular integrals and differentiability properties of functions.} Princeton University Press   (1970).

\bibitem{strichartz94}  R. Strichartz,  \textit{A guide to distribution theory and Fourier transforms.} Studies in Advanced Mathematics. CRC Press, Boca Raton, FL, (1994).

\bibitem{taylor96}
M. Taylor,
\textit{Partial differential equations. I. }
Basic theory. Applied Mathematical Sciences, 115.
Springer-Verlag, New York, 1996.

 \end{thebibliography}
\end{document}